\documentclass[10pt]{amsart}
\usepackage{latexsym,enumerate}
\usepackage{amsmath,amsthm,amsopn,amstext,amscd,amsfonts,amssymb}
\usepackage{color}
\usepackage{fancybox}
\usepackage{epsfig}
\usepackage{float}
\usepackage{ulem} %para poder tachar: \sout{loquequierotachar}

\usepackage[active]{srcltx}

\usepackage{color}

\newcommand{\R}{\mathbb{R}}
\newcommand{\N}{\mathbb{N}}

\newcommand{\weakly}{\rightharpoonup}

\newcommand{\car}{{\raise4pt\hbox{$\chi$}}}

\newcommand{\sg}{{\rm \; sign \;}}

\newcommand{\sop}{{\rm supp\,}}

\newcommand{\Div}{\hbox{\rm div\,}}

\newcommand{\qin}{\qquad\mbox{ in }\quad}
\newcommand{\qon}{\qquad\mbox{ on }\quad}

\newcommand{\dis }{{\mathcal D}' }
\newcommand{\z }{{\bf z}}
\newcommand{\DM }{\mathcal{DM}^\infty }
\newenvironment{pf}{\noindent{\sc Proof}.\enspace}{\rule{2mm}{2mm}\medskip}

\newtheorem{Theorem}{Theorem}[section]
\newtheorem{Corollary}[Theorem]{Corollary}
\newtheorem{Definition}[Theorem]{Definition}

\newtheorem{Proposition}[Theorem]{Proposition}

\newtheorem{remark}[Theorem]{Remark}
\newtheorem{Example}[Theorem]{Example}
\newcommand{\res}{\!\!\mathop{\hbox{
                                \vrule height 7pt width .5pt depth 0pt
                                \vrule height .5pt width 6pt depth 0pt}}
                                \nolimits}
\newcommand{\norma}[2]{\|#1\|_{\lower 4pt \hbox{$\scriptstyle #2$}}}

\usepackage{multirow}
\usepackage{array}

\begin{document}

\title[The 1--Laplacian equation with a singular lower order term]{Elliptic problems involving the 1--Laplacian\\ and a singular lower order term}

\author[V. De Cicco, D. Giachetti and  S. Segura de Le\'on]
{V. De Cicco, D. Giachetti and  S. Segura de Le\'on}

\address{V. De Cicco: Dipartimento di Scienze di Base e Applicate per l'Ingegneria,
Sapienza--Universit\`a di Roma,
Via A. Scarpa 16, 00161 Roma, Italy.\newline
{\it E-mail address:} {\tt virginia.decicco@sbai.uniroma1.it }}

\address{D. Giachetti: Dipartimento di Scienze di Base e Applicate per l'Ingegneria,
Sapienza--Universit\`a di Roma,
Via A. Scarpa 16, 00161 Roma, Italy.\newline
{\it E-mail address:} {\tt daniela.giachetti@sbai.uniroma1.it }}

\address{S. Segura de Le\'on: Departament d'An\`{a}lisi Matem\`atica,
Universitat de Val\`encia,\hskip1cm
Dr. Moliner 50, 46100 Burjassot, Valencia, Spain.\newline
{\it E-mail address:}  {\tt sergio.segura@uv.es }}

\thanks{}

%\date{\empty}

\begin{abstract}
%{{\color{red}
  This paper is concerned with the Dirichlet problem for an equation involving the $1$--Laplacian operator $\Delta_1 u:=\Div\left(\frac{Du}{|Du|}\right)$ and having a singular term of the type $\frac{f(x)}{u^\gamma}$. Here $f\in L^N(\Omega)$ is nonnegative, $0<\gamma\le1$ and $\Omega$ is a bounded domain with Lipschitz--continuous boundary. We prove an existence result for a concept of solution conveniently defined. The solution is obtained as limit of solutions of $p$--Laplacian type problems. Moreover, when $f(x)>0$ a.e., the solution satisfies those features that might be expected as well as a uniqueness result. We also give explicit 1--dimensional examples that show that, in general, uniqueness does not hold. We remark that the Anzellotti theory of $L^\infty$--divergence--measure vector fields must be extended to deal with this equation.
%}}
\end{abstract}

\maketitle

%\tableofcontents

\section{Introduction}
In the present paper we deal with the Dirichlet problem for
equations involving the $1$--Laplacian and a singular lower order term. This equation, at least in the case $f>0$ a.e. in $\Omega$, looks like
\begin{equation}\label{prob-prin}
  \left\{\begin{array}{ll}
  \displaystyle -\Div\Big(\frac{Du}{|Du|}\Big)=\frac{f(x)}{u^\gamma}&\hbox{in
   }\Omega\,,\\[3mm]
   u=0 &\hbox{on }\partial\Omega\,,
  \end{array}\right.
\end{equation}
where $\Omega\subset\R^N$ is a bounded open set with Lipschitz
boundary $\partial\Omega$, $0<\gamma\le1$ and $f$ is a function belonging to $L^{N}(\Omega)$. Actually we are interested to the case of $f\geq 0$. The natural space to study problems where the $1$--Laplacian appears is $BV(\Omega)$, the space of functions $u\in L^1(\Omega)$ whose distributional derivatives are Radon measures with finite total variation.

We point out that problems involving a lower order singular term like that appearing in problem \eqref{prob-prin} have already been studied in the Laplacian or the $p$--Laplacian setting.
 There is a wide literature dealing with problem
 \begin{equation}\label{eqprimai}
\begin{cases}
\displaystyle -\Delta_p u =\frac{f(x)}{u^\gamma}  & \mbox{in }  \Omega,\\
u = 0 & \mbox{on} \; \partial \Omega\\
\end{cases}
\end{equation}
for $p>1$ and $f\geq 0$.
%We will not pretend to give here a complete list of references and we will concentrate on some papers which seem to us to be  the most significant, also refering the interested reader to the references quoted there.
%
The problem \eqref{eqprimai} for $p=2$ was initially proposed in 1960 in the pioneering work \cite{FM} by Fulks and Maybee as a model for several physical situations. This problem was then studied by many authors, among which we will quote the works of  Stuart \cite{St},  Crandall, Rabinowitz and Tartar \cite{BCR},  Lazer and McKenna \cite{LM},  Boccardo and Orsina  \cite{BO}, Coclite and Coclite \cite{CoCo},
Arcoya and Moreno--M\'erida \cite{AMo},  Oliva and Petitta \cite{OlPe}, and Giachetti, Martinez--Aparicio and Murat in \cite{GMM1}, \cite{GMM2}, \cite{GMM} and \cite{GMM2bis}.
% In most of the above quoted papers the authors look for a strong solution and use sub- and super-solutions.  In particular in \cite{LM}  A.C.~Lazer and P.J.~McKenna work in $C^{2,\alpha}(\Omega)$ and $W^{2,q}(\Omega)$  and use methods of sub- and super-solutions, proving  that when $\displaystyle F(x,s)=\frac{f(x)}{s^\gamma}$ with $\gamma>0$, $f\in C^\alpha(\overline{\Omega})$ and $f(x)\geq f_0>0$ in a $C^{2,\alpha}$ domain $\Omega$, then one has
%$
%c_1\phi_1(x)\leq u(x)^{\frac{\gamma+1}{2}}\leq c_2\phi_1(x)
%$
%for two constants $0<c_1<c_2$, where $\phi_1$ is the first (positive) eigenfunction of $-div\, A(x)D$ in $H_0^1(\Omega)$;
%M.G.~Crandall, P.H.~Rabinowitz and L.~Tartar \cite{BCR} study the behaviours of $u(x)$ and $|Du(x)|$ at the boundary. Let us finally note that C.~Stuart \cite{St},  as well as  M.G.~Crandall, P.H.~Rabinowitz and L.~Tartar \cite{BCR}, do not assume that $F(x,s)$ is nonincreasing in $s$.

In particular, in \cite{BO} the authors studied the problem in the framework of weak solutions in the sense of distributions for $f$ belonging to Lebesgue spaces, or to the space of Radon measures, and they prove existence and regularity as well as non existence results.
On the other hand, in \cite{GMM1}, \cite{GMM2} and \cite{GMM}  the authors deal with fairly general  singular problems and they find weak solutions belonging, for $\gamma\leq1$, to the energy space.

% In this work the strong maximum principle and the nonincreasing character of the function $\displaystyle F(x,s)=\frac{f(x)}{s^\gamma}$ with respect to $s$ play prominent roles. The solution $u$ to problem \eqref{eqprimai} is indeed required to satisfy $u(x)\geq c(\omega,u)>0$ on every open set $\omega$ with $\overline{\omega}\subset \Omega$, and the framework in which the solution is searched for is the Sobolev's space  $H_0^1(\Omega)$ or $H^1_{\mbox{\tiny loc}}(\Omega)$. Then L.~Boccardo and J.~Casado-D\'iaz  proved in \cite{BC} the uniqueness of the solutions obtained by approximation and the stability of the solution with respect to the $G$-convergence of a sequence of matrices $A^\eps(x)$ which are equicoercive and equibounded.

For a variational approach to the problem \eqref{eqprimai} in the case $p>1$, see for instance Canino and Degiovanni \cite{CDG}.
%Arcoya, Boccardo and Orsina \cite{BO}.

Due to the methods we are going to use in our paper (approximation of our problem with problems driven by the $p$--Laplacian), we quote now
some papers dealing with existence results for problem \eqref{eqprimai} in the case $p>1$. We refer to Giacomoni, Schindler and Takac  \cite{GST}, Perera and Silva \cite{PE},  Mohammed \cite{Mo}, Loc Hoang and Schmitt \cite{LHS}, De~Cave \cite{DCA}, Canino, Sciunzi and Trombetta \cite{CST} and Mi \cite{Mi}.

Early papers devoted to the Dirichlet problem for equations involving the $1$--Laplacian operator include \cite{K, HI, ABCM, BCN, D, CT}. The interest in this setting comes out  from an optimal design problem in the theory of torsion and from the level set formulation of the Inverse Mean Curvature Flow.  On the other hand, it also appears in the variational approach to image restoration. Indeed, total variation minimizing models have become one of the most popular and successful methodology for image restoration since the introduction of the ROF model by Rudin, Osher and Fatemi in \cite{ROF}. In this paper a variational problem involving the total variation operator is considered. It was designed with the explicit goal of preserving sharp discontinuities (edges) in images while removing noise and other unwanted fine scale detail. The $1$--Laplacian operator emerges through the subdifferential of the total variation. Since that paper, there has been a burst in the application of the total variation regularization to many different image processing problems which include inpainting, blind deconvolution or multichannel image segmentation (for a review on this topic we refer to \cite{CEPY}, see also \cite{ACM}).

To deal with the $1$--Laplacian, a first difficult occurs in defining the quotient $\displaystyle \frac{Du}{|Du|}$, being $Du$ a Radon measure. It was overcome by Andreu, Ballester, Caselles and Maz\'on in \cite{ABCM} through the Anzellotti theory of pairings of $L^\infty$--divergence--measure vector fields and the gradient of a BV--function (see \cite{ACM}). In their definition appears a vector field $\z\in L^\infty(\Omega;\R^N)$ such that $\|\z\|_\infty\le 1$ and $(\z, Du)=|Du|$, so that $\z$ plays the role of the above rate.

Since the Dirichlet boundary condition does not hold in the usual trace form, they also introduce in \cite{ABCM} a weak formulation: $[\z, \nu]\in\sg(-u)$, where $[\z,\nu]$ stands for the weak trace on $\partial\Omega$ of the normal component of $\z$.

It is worth observing that the related problem having no singular lower order term
\begin{equation}\label{prob-secon}
  \left\{\begin{array}{ll}
  \displaystyle -\Div\Big(\frac{Du}{|Du|}\Big)=f(x)&\hbox{in
   }\Omega\,,\\[3mm]
   u=0 &\hbox{on }\partial\Omega\,,
  \end{array}\right.
\end{equation}
with $f\in L^N(\Omega)$,
has particular features and there exists a BV--solution only if $\|f\|_N$ is small enough (see \cite{CT, MST1}).

Our main results show that the term $\frac{f(x)}{u^\gamma}$ has a regularizing effect. Indeed, in Theorem \ref{teoexist} below, we prove that there exists a bounded solution to problem \eqref{prob-prin}, in the sense of Definition \ref{def1}, for every nonnegative datum $f\in L^N(\Omega)$.

Moreover,  in Theorem \ref{teoexist1}, we see the improved features satisfied by the solution when $f$ is strictly positive, which lead to uniqueness (Theorem \ref{uniq}). In summary, the bigger the data, better features has the solution.

Being the source term singular on the set $\{u=0\}$, this set is crucial. If $p=1,$ we are able to prove (as in the $p>1$ case) that \begin{equation*}
  \{x\in\Omega: u(x)=0\}\subset \{x\in\Omega: f(x)=0\}\,,
\end{equation*}
except for a set of zero Lebesgue measure and we will prove that this set has locally finite perimeter. The conditions satisfied by solutions for non strictly positive data include
\begin{enumerate}[$(a)$]
\item $\displaystyle \frac f{u^\gamma} \in L^1_{\rm{loc}}(\Omega)$,

\item  $\displaystyle-(\Div \z)\chi^*_{\{u>0\}}  =\frac f{u^\gamma}$ in $ \dis (\Omega)$,
\end{enumerate}
where $\z$ is the vector field defined above and $\chi^*_{\{u>0\}}$ is the precise representative (in the sense of the $BV$--function) of the characteristic function $\chi_{\{u>0\}}$ of the set ${\{u>0\}}$. We will prove that ${\{u>0\}}$ is a set of locally finite perimeter and so $\chi_{\{u>0\}}$ is a locally $BV$--function.
Note that,
by using that
\begin{equation*}
\Div\big(\z\chi_{\{u>0\}}\big)=\big(\Div \z\big)\chi^*_{\{u>0\}}+(\z,D \chi_{\{u>0\}}),
\end{equation*} the equation (b)
%$$\displaystyle-(div \,z)\chi^*_{\{u>0\}}  =\frac {f(x)}{u^\gamma}$$
can be written as
\begin{equation*}
-\Div\big(\z\chi_{\{u>0\}}\big)+|D \chi_{\{u>0\}}|=\frac {f(x)}{u^\gamma},
\end{equation*}
where the left hand side is a sum of an operator in divergence form and an additional term $|D \chi_{\{u>0\}}|$, which is a measure concentrated on the reduced boundary $\partial^*\{u>0\}$ (or equivalently on the reduced boundary $\partial^*\{u=0\}$).

Since we do not know that the measure $\Div\z$ has finite total variation, we cannot directly apply Anzellotti's theory. As a consequence, one of our tasks in analyzing problem \eqref{prob-prin} will be slightly extend it.

The proof of our results is obtained passing to the limit in singular problems with principal part the $p$--Laplacian, $p>1$ as $p$ tends to $1.$ To this aim, we need preliminarly to prove that the corresponding singular approximating problems admit a weak solution (see Theorem \ref{pproblem}).

Next we will show that the family $(u_p)_{p>1}$ is bounded in the $BV$ norm, hence there exists a function $u \in BV(\Omega)$ such that,
up to a
subsequence,
\begin{equation*}
u_p \to u\ {\text{ in }}\  L^1(\Omega)
\end{equation*}
and
\begin{equation*}
\nabla u_p \rightharpoonup Du\ {\text{ weakly} }^* {\text{\ as measures}}.
\end{equation*}
Similarly, for the vector field $\z$, we have that
$(|\nabla u_p|^{p-2}\nabla u_p)_p$ is bounded in $L^q(\Omega;\R^N)$ for any $1\leq q<\infty.$ Hence there exists a vector field $\z$ satisfying
\begin{equation*}
|\nabla u_p|^{p-2}\nabla u_p\rightharpoonup \z\,,\quad\hbox{weakly in }L^q(\Omega;\R^N)\,,\quad\forall 1\leq q<\infty\,.
\end{equation*}
Such functions $u$ and $\z$ are those satisfying  (a) and (b) as well as the other conditions of Definition \ref{def1}.

Let us remark that if $f\equiv0$ the solutions are $1$--harmonic functions; if the boundary datum is continuous and the domain satisfies some additional geometrical assumptions, the solution is unique (see Sternberg, Williams and Ziemer in \cite{SWZ}) and it is the limit of $p$--harmonic functions (see Juutinen in \cite{J}). 
%We find these known results as a particular case.
If $f$ is not strictly positive and $f$ not identically zero there are 1--dimensional examples of non uniqueness (see Section 6).

The plan of this paper is the following. Section 2 is dedicated to
preliminaries, we introduce our notation, some properties of
the space $BV(\Omega)$ and an extension of the Anzellotti theory, including a Green formula. The starting point of our main result is studied in Section 3, which is concerned with solutions to problem \eqref{eqprimai}.  Section 4 is devoted to the study of existence for general nonnegative data, while Section 5 deals with strictly positive data. Finally, Section 6 provides some explicit 1--dimensional solutions that show that, in general, uniqueness does not hold.

\section*{Acknowledgements}
This research has been partially supported  by the Spanish Mi\-nis\-te\-rio de Eco\-no\-m\'{\i}a y Competitividad and
FEDER, under project MTM2015-70227-P. This paper was started during a stay of the third author in Sapienza--Universit\`a di Roma supported with a grant as visiting professor by Istituto Nazionale di Alta Matematica Francesco Severi.
%he also thanks Dipartimento di Scienze di Base e Applicate per l'Ingegneria for its warm hospitality. 
The authors also thank Francesco Petitta for some useful suggestions.

\section{Preliminaries}

In this Section we will introduce some notation and auxiliary results which will be used
throughout this paper. In what follows, we will consider $N\ge2$,
and $\mathcal H^{N-1}(E)$ will denote the $(N - 1)$--dimensional
Hausdorff measure of a set $E$ and $|E|$ its
Lebesgue measure.

  In this paper, $\Omega$ will always denote an open bounded subset of
  $\R^N$ with Lipschitz boundary. Thus, an outward normal unit
  vector $\nu(x)$ is defined for $\mathcal H^{N-1}$--almost every
  $x\in\partial\Omega$.
 We will make use of the usual Lebesgue and Sobolev
 spaces, denoted by $L^q(\Omega)$  and $W_0^{1,p}(\Omega)$,
 respectively. On the other hand, Lebesgue spaces with respect to a Radon measure $\mu$ are denoted by $L^q(\Omega, \mu)$.

We recall that for a Radon measure $\mu$ in $\Omega$ and a Borel set
$A\subseteq\Omega$ the measure $\mu\res A$ is defined by $(\mu\res A)(B)=\mu(A\cap B)$
for any Borel set $B\subseteq\Omega$. If a measure $\mu$ is such that $\mu
 = \mu \res A$ for a certain Borel set $A$, the measure $\mu$ is
 said to be concentrated on $A$.

The truncation function will be use throughout this paper. Given
$k>0$, it is defined by
\begin{equation*}\label{trun}
    T_k(s)=\min\{|s|, k\}\sg (s)\,,
\end{equation*} for all $s\in\R$. Moreover we will denote by $G_k(s)$ the function defined by
$$
G_k(s)=s-T_k(s)
$$
for all $s\in\R$.
\subsection{The energy space}

The space of all functions of finite variation, that is the space
of those $u\in L^1(\Omega)$ whose distributional gradient is a Radon
measure with finite total variation, is denoted by $BV(\Omega)$.
This is the natural energy space to study the problems we are interested
in. It
is endowed with the norm defined by
 $$ \|u\|=\int_\Omega |u|\, dx+ \int_\Omega|Du|\,,$$
 for any $u\in BV(\Omega)$.
 We recall that the notion of trace can be extended to any $u\in BV(\Omega)$ and this fact allows us to interpret it as the boundary values of $u$ and to write $u\big|_{\partial \Omega}$. Moreover, the trace defines a linear bounded operator $BV(\Omega)\to  L^1(\partial\Omega)$ which is onto. Using the trace, we have available an equivalent norm, which we will use in the sequel. It is given by
 $$\displaystyle \|u\|_{BV(\Omega)}=\int_{\partial\Omega}
|u|\, d\mathcal H^{N-1}+ \int_\Omega|Du|\,.$$

For every
$u \in BV(\Omega)$, the Radon measure $Du$ is decomposed into its
absolutely continuous and singular parts with respect to the
Lebesgue measure: $Du = D^a u + D^s u$.
We denote by $S_u$ the set of all $x\in\Omega$ such that $x$ is
not a Lebesgue point of $u$, that is, $x\in\Omega\backslash S_u$ if there exists $\tilde{u}(x)$ such that
$$\lim_{\rho \downarrow 0} \frac{1}{|B_{\rho}(x)|}
\int_{B_{\rho}(x)} \vert u(y) - \tilde{u}(x) \vert \, dy = 0\,.$$
 We say that $x \in \Omega$ is an {\it
approximate jump point of } $u$ if there exist two real numbers $u^+(x)
>u^-(x) $ and $\nu_u(x) \in S^{N-1}$ such that
$$\lim_{\rho \downarrow 0} \frac{1}{|B_{\rho}^+(x,\nu_u(x))|}
\int_{B_{\rho}^+(x,\nu_u(x))} \vert u(y) - u^+(x) \vert \, dy = 0\,,$$
$$\lim_{\rho \downarrow 0} \frac{1}{|B_{\rho}^-(x,\nu_u(x))|}
\int_{B_{\rho}^-(x,\nu_u(x))} \vert u(y) - u^-(x) \vert \, dy = 0\,,$$ where
$$B_{\rho}^+(x,\nu_u(x)) = \{ y \in B_{\rho}(x) \ : \ \langle y - x, \nu_u(x) \rangle >
0 \} $$ and
$$B_{\rho}^-(x,\nu_u(x)) = \{ y \in B_{\rho}(x) \ : \ \langle y - x, \nu_u(x) \rangle <
0 \}\,.$$
 We denote by $J_u$ the set of all approximate jump points of
$u$. By the Federer--Vol'pert Theorem \cite[Theorem 3.78]{AFP}, we know that $S_u$ is countably $\mathcal H^{N-1}$--rectifiable and $\mathcal H^{N-1}(S_u \backslash J_u) = 0$.
%{\color{red}{
%Moreover, $Du \res J_u = (u^+ - u^-) \nu_u \mathcal H^{N-1} \res J_u$. Using
%$S_u$ and $J_u$, we may split $D^su$ in two parts: the {\it jump}
%part $D^j u$ and the {\it Cantor} part $D^c u$ defined by
%$$D^ju = D^su \res J_u  \ \ \ {\rm and} \ \ D^c u = D^su \res (\Omega \backslash S_u)\,.$$
%Then, we have
%$$D^j u = (u^+ - u^-) \nu_u \mathcal H^{N-1} \res J_u\,.$$
%Moreover, if $x \in J_u$, then $\nu_u(x) = \frac{Du}{| D u |}(x)$ and
%$\frac{Du}{| D u |}$ is the Radon--Nikod\'ym derivative of $Du$
%with respect to its total variation $| D u |$.
%}}

The precise representative
$u^* : \Omega \backslash(S_u \backslash J_u) \rightarrow \R$ of
$u$ is defined as equal to $\tilde{u}$ on $\Omega \backslash S_u$
and equal to $\frac{u^- + u^+}{2}$ on $J_u$. It is well known (see
for instance \cite[Corollary 3.80]{AFP}) that if $\rho$ is a
symmetric mollifier, then the mollified functions $u \star
\rho_{\epsilon}$ pointwise converges to $u^*$ in its domain.

  A compactness result in $BV(\Omega)$ will be used in what follows.
  It states that every sequence that is bounded in $BV(\Omega)$ has a
subsequence which strongly converges in $L^1(\Omega)$ to a certain
$u\in BV(\Omega)$ and the subsequence of gradients $*$--weakly
converges to $Du$ in the sense of measures.

 To pass to the limit we will often apply that some functionals defined on $BV(\Omega)$ are
 lower semicontinuous with respect to the convergence in $L^1(\Omega)$.
  The most important are
  the functionals
  defined by
  \begin{equation*}\label{semc}
 u\mapsto\int_{\Omega}|Du|
 \end{equation*}
 and
   \begin{equation*}\label{semcon}
  u\mapsto\int_\Omega|Du|+\int_{\partial\Omega}|u|\,d\mathcal
  H^{N-1}\,.
  \end{equation*}
  In the same way,
  it yields that each $\varphi\in C_0^1(\Omega)$ with
  $\varphi\ge0$ defines a functional
  $$
  u\mapsto\int_{\Omega}\varphi\,|Du|\,,
  $$
  which is lower semicontinuous in $L^1(\Omega)$.

For further information on functions of bounded variation, we refer to \cite{AFP, EG, Zi}.

\subsection{A generalized Green formula}
The theory of $L^\infty$--divergence--measure vector fields is due
to Anzellotti \cite{An} and to Chen and
Frid \cite{CF}. In spite of their different points of view,
both approaches introduce, under some hypotheses, the ``dot product" of a bounded vector field $\z$, whose divergence is a Radon measure, and the gradient $Dv$ of $v\in BV(\Omega)$ through a pairing $(\z,Dv)$ which defines a Radon measure. However, they differ in handling this concept. They also define the normal trace of a vector field through the boundary and establish a generalized
Gauss--Green formula.

From now on, we denote by $\DM(\Omega)$ the space of all
 vector fields $\z\in L^\infty(\Omega;\R^N)$ whose divergence
 in the sense of distribution is a Radon measure with finite total variation, i.e., $\z\in \DM(\Omega)$ if and only if
 $\Div \z$ is a finite Radon measure belonging to $W^{-1,\infty}(\Omega)$. Moreover, $\DM_{\rm{loc}}(\Omega)$ stands for those vector fields $\z\in L^\infty(\Omega;\R^N)$ which belong to $\DM(\omega)$ for every open $\omega\subset\!\subset\Omega$.

To define $(\z,Dv)$ in Anzellotti's theory, we need some compatibility conditions, such as $\Div\z\in L^1(\Omega)$ and $v\in BV(\Omega)\cap L^\infty(\Omega)$, or $\Div\z$ a Radon measure with finite total variation and $v\in BV(\Omega)\cap L^\infty(\Omega)\cap C(\Omega)$. Anzellotti's definition of $(\z, Dv)$ can be extended to the case which $\Div\z$ is a Radon measure with finite total variation and $v\in BV(\Omega)\cap L^\infty(\Omega)$ (see \cite[Appendix A]{MST2} and \cite[Section 5]{C}). We are following these papers for a slightly further extension to the spaces $\DM_{\rm{loc}}(\Omega)$ and $BV_{\rm{loc}}(\Omega)$. The starting point is also the following result proved in \cite{CF}.

   \begin{Proposition}\label{absolcont}
  For every $\z\in\DM(\Omega)$,
  the measure $\mu=\hbox{\rm div\,}\z$ is absolutely continuous with respect to
 $\mathcal H^{N-1}$, that is,  $|\mu|\ll\mathcal H^{N-1}$.
 \end{Proposition}

 Consider now $\mu=\hbox{\rm div}\,\z$ with $\z\in\DM_{\rm{loc}}(\Omega)$ and let $v\in
BV_{\rm{loc}}(\Omega)$. Since the precise representative $v^*$ is equal
$\mathcal H^{N-1}$--a.e. to the Borel function $\lim_{\epsilon\to0}v\star\rho_\epsilon$, then one deduces from
Proposition \ref{absolcont}, that $v^*$ is equal
$\mu$--a.e. to a Borel function so that the precise representative of every BV--function is $\mu$--measurable. Assume that $v^*\in L^1(\Omega, \mu)$ and
 define a distribution by the following expression
 \begin{equation}\label{dist1}
\langle(\z,Dv),\varphi\rangle:=-\int_\Omega v^*\varphi\,d\mu-\int_\Omega
 v\z\cdot\nabla\varphi \,dx,\quad \varphi\in C_0^\infty(\Omega)\,.
\end{equation}
 Every term is well defined since
  $v\in BV_{\rm{loc}}(\Omega)\cap L^1(\Omega,\mu)$
  and
 $\z \in L^{\infty}(\Omega, \R^N)$. We point out that the definition of $(\z,Dv)$ depends on the precise representative of $v$; so that if we choose another representative, the distribution will become different.

We next see that this distribution
is actually a Radon measure having locally finite total variation. The proofs are similar to those in \cite{MST2} or \cite{C}.

\begin{Proposition}\label{prop}
  Let $v\in BV_{\rm{loc}}(\Omega)\cap L^1(\Omega,\mu)$ and
 $\z \in \DM_{\rm{loc}}(\Omega)$.
 Then the distribution $(\z, Dv)$ defined previously satisfies
  \begin{equation}\label{ec:2}
  |\langle   (\z, Dv), \varphi\rangle| \le \|\varphi\|_\infty \| \z
\|_{L^{\infty}(U)} \int_{U} |Dv|
  \end{equation}
for all open set $U \subset\!\subset \Omega$ and for all $\varphi\in
C_0^\infty(U)$.
 \end{Proposition}

   \begin{pf}
  If $U\subset \Omega$ is an open set and $\varphi\in
 C_0^\infty(U)$, then it was proved in \cite{MST2} that
  \begin{equation}\label{ec:3}
      |\langle(\z, DT_k(v)),\varphi\rangle|\le\|\varphi\|_\infty
  \|\z \|_{L^\infty(U)}\int_U|DT_k(v)|\le\|\varphi\|_\infty
  \|\z \|_{L^\infty(U)}\int_U|Dv|
\end{equation}
 holds for every $k>0$.
 On the other hand,
 $$
 \langle(\z, DT_k(v)),\varphi\rangle=-\int_\Omega T_k(v)^*\varphi
 \, d\mu-\int_\Omega T_k(v)\z\cdot\nabla\varphi\, dx\,.
 $$
 We may let $k\to\infty$ in each term on the right hand side, due
 to $v^*\in L^1(\Omega, \mu)$ and $v\in L^1(\Omega)$.
 Therefore,
 $$
 \lim_{k\to\infty}\langle(\z, DT_k(v)),\varphi\rangle=\langle(\z,
 Dv),\varphi\rangle\,,
 $$
 and so \eqref{ec:3} implies \eqref{ec:2}.
 \end{pf}

 \begin{Corollary}
  The distribution $(\z, Dv)$ is a Radon measure. It and its total variation $\vert (\z, Dv) \vert$ are absolutely
continuous with respect to the measure $\vert Dv \vert$ and
$$\left\vert \int_{B}  (\z, Dv) \right\vert \leq \int_{B} \vert (\z, Dv) \vert \leq \Vert \z
\Vert_{L^{\infty}(U)} \int_{B} \vert Dv \vert
$$
holds for all Borel sets $B$ and for all open sets $U$ such that
$B \subset U \subset \Omega$.

In particular, if $v\in BV(\Omega)$, then the measure $(\z, Dv)$ has finite total variation.
 \end{Corollary}

 Moreover, going back to \eqref{dist1}, we can conclude that the following proposition holds.

 \begin{Proposition}\label{nuova}
Let  $\z\in\DM_{\rm{loc}}(\Omega)$ and let $v\in
BV_{\rm{loc}}(\Omega)\cap L^\infty(\Omega)$. Then $\z v\in\DM_{\rm{loc}}(\Omega)$. Moreover the following formula holds in the sense of measures
\begin{equation}
\label{ALPHA}
\Div(\z v)=(\Div \z) v^*+(\z,Dv).
\end{equation}
 \end{Proposition}

 On the other hand, for every $\z \in \mathcal{DM}^{\infty}(\Omega)$,
 a weak trace on $\partial \Omega$ of the
normal component of  $\z$ is defined in \cite{An} and denoted by
$[\z, \nu]$. Moreover, it satisfies
\begin{equation}\label{des1}
  \|[\z,\nu]\|_{L^\infty(\partial\Omega)}\le \|\z\|_{L^\infty(\Omega)}\,.
\end{equation}
We explicitly point out that if $\z \in \mathcal{DM}^{\infty}(\Omega)$ and $v\in BV(\Omega)\cap L^\infty(\Omega)$, then
\begin{equation}\label{des2}
v[\z,\nu]=[v\z,\nu]
\end{equation}
 holds (see \cite[Lemma 5.6]{C} or \cite[Proposition 2]{ADS}). Having the ``dot product" $(\z , Dv)$ and the normal component $[\z, \nu]$, we may already prove our first Green formula.

\begin{Proposition}
Let $\z \in \mathcal{DM}^{\infty}(\Omega)$, $v\in BV(\Omega)$ and assume that $v^*\in L^1(\Omega,\mu)$.
 With the above definitions, the
 following Green formula holds
\begin{equation}\label{GreenI}
\int_{\Omega} v^* \, d\mu + \int_{\Omega} (\z, Dv) =
\int_{\partial \Omega} [\z, \nu] v \ d\mathcal H^{N-1}\,.
\end{equation}
 \end{Proposition}

 \begin{pf}
Applying the Green formula proved in \cite[Theorem A.1]{MST2} or in \cite[Theorem 5.3]{C}, we obtain
  \begin{equation}\label{ec:4}
\int_{\Omega} T_k(v)^* \, d\mu + \int_{\Omega} (\z, DT_k(v)) =
\int_{\partial \Omega} [\z, \nu] T_k(v) \ d\mathcal H^{N-1}\,,
\end{equation}
for every $k>0$.  Note that the same argument appearing in the proof of Proposition \ref{prop} leads to
$$
\lim_{k\to\infty}\int_{\Omega} (\z, DT_k(v)) =\int_{\Omega} (\z,
Dv)\,.
$$
We may take limits in the other terms since $v^*\in L^1(\Omega,
\mu)$ and $v\in L^1(\partial\Omega)$. Hence, letting $k$ go to
$\infty$ in \eqref{ec:4}, we get \eqref{GreenI}.
\end{pf}

Observe that we may apply \eqref{GreenI} to a vector field $\z \in \mathcal{DM}^{\infty}(\Omega)$ and the constant $v\equiv 1$. Since $(\z, Dv)=0$, we obtain
\begin{equation}\label{GreenII}
\int_{\Omega}  \Div\z =
\int_{\partial \Omega} [\z, \nu] \ d\mathcal H^{N-1}\,.
\end{equation}

This fact and \eqref{ALPHA} are enough to prove the Green formula we will apply in what follows.

\begin{Proposition}\label{poiu}
Let $\z \in \mathcal{DM}_{\rm{loc}}^{\infty}(\Omega)$ and set $\mu=\Div\z$. Let $v\in BV(\Omega)\cap L^\infty(\Omega)$ be such that that $v^*\in L^1(\Omega,\mu)$.
Then $v\z\in \mathcal{DM}^{\infty}(\Omega)$ and the following Green formula holds
\begin{equation}\label{GreenIII}
\int_{\Omega} v^* \, d\mu + \int_{\Omega} (\z, Dv) =
\int_{\partial \Omega} [v\z, \nu] \ d\mathcal H^{N-1}\,.
\end{equation}
 \end{Proposition}

 \begin{pf}
Taking into account that $v^*\in L^1(\Omega,\mu)$ implies that
\[
\Div(v\z)=v^*\Div\z+(\z,Dv)
\]
is a Radon measure with finite variation, and so $v\z\in \mathcal{DM}^{\infty}(\Omega)$ (observe that we assume $v\in L^\infty(\Omega)$).
It follows now from \eqref{GreenII} that \eqref{GreenIII} holds.
\end{pf}

%{{\color{red}
We will deal with $\z\in \DM_{\rm{loc}}(\Omega)$ such that the product $v\z\in \DM(\Omega)$ for some $v\in BV(\Omega)\cap L^\infty(\Omega)$.
To infer the Dirichlet boundary condition, we will use the following result.

\begin{Proposition}
Let $\z \in \mathcal{DM}_{\rm{loc}}^{\infty}(\Omega)$ and $v\in BV(\Omega)\cap L^\infty(\Omega)$. If $v\z\in \DM(\Omega)$, then
\begin{equation}\label{des3}
  |[v\z,\nu]|\le |v\big|_{\partial\Omega}|\,\|\z\|_{L^\infty(\Omega)}\,\quad\mathcal H^{N-1}\hbox{--a.e. on }\partial\Omega\,.
\end{equation}
 \end{Proposition}

 \begin{pf}
 We first claim that,
 \begin{equation}\label{claim}
   \|\phi[v\z,\nu]\|_{L^\infty(\partial\Omega)}\le \|\phi v\big|_{\partial\Omega}\|_{L^\infty(\partial\Omega)}\,\|\z\|_{L^\infty(\Omega)}
 \end{equation}
 holds for every nonnegative $\phi\in L^\infty(\partial\Omega)$.

 To this end, let $\varphi\in BV(\Omega)\cap L^\infty(\Omega)$ satisfy $\varphi\ge0$ and $\varphi\big|_{\partial\Omega}=\phi$.
 Applying \eqref{des1} and \eqref{des2}, we may manipulate as follows
 \begin{equation*}
   \|\phi[v\z,\nu]\|_{L^\infty(\partial\Omega)}=\|[\varphi v\z,\nu]\|_{L^\infty(\partial\Omega)}\le \|\varphi v \z\|_{L^\infty(\Omega)}
   \le \|\varphi v\|_{L^\infty(\Omega)}\| \z\|_{L^\infty(\Omega)}\,.
 \end{equation*}
 By \cite[Lemma 5.5]{An}, we find a sequence $(\varphi_n)_n$ in
 $W^{1, 1}(\Omega)\cap C(\Omega)\cap L^\infty(\Omega)$ such that
 \begin{equation*}\begin{array}{ll}
 \hbox{\rm (1) }
 \varphi_n|_{\partial\Omega}=\phi\,.\\
\\
 \hbox{\rm (2) } \|\varphi_n\|_{L^\infty(\Omega)}=\|\phi\|_{L^\infty(\partial\Omega)}\,.\\
\\
 \hbox{\rm (3) } \varphi_n(x)=0,  \qquad\hbox{ if } dist(x,\partial\Omega)>\frac1n\,,
\hskip5cm
 \end{array}\end{equation*}
 for all $n\in\N$.
   It follows from $\|\phi[v\z,\nu]\|_{L^\infty(\partial\Omega)}\le \|\varphi_n v\|_{L^\infty(\Omega)}\| \z\|_{L^\infty(\Omega)}$ for all $n\in\N$ that \eqref{claim} holds true.

   To prove \eqref{des3}, assume to get a contradiction that there exist $\epsilon>0$ and a measurable set $E\subset\partial\Omega$ such that $\mathcal H^{N-1}(E)>0$ and
   \[
   |[v\z,\nu]|\ge |v\big|_{\partial\Omega}|\,\|\z\|_{L^\infty(\Omega)}+\epsilon\,,\quad\hbox{on }E\,.
   \]
   Letting $\phi=\chi_E$, we deduce
   \[
   \|\phi[v\z,\nu]\|_{L^\infty(\partial\Omega)}\ge \|\phi v\big|_{\partial\Omega} \|_{L^\infty(\partial\Omega)}\|\z\|_{L^\infty(\Omega)}+\epsilon\,,
   \]
   which contradicts \eqref{claim}.
 \end{pf}
%}}
\section{Weak solution for $p$--Laplacian type problems}
For every $p>1$ let us consider the following problem
\begin{equation}
    \left \{
        \begin{array}{cl}
            \displaystyle-\Delta_p \left(u\right) = \frac{f(x)}{u^\gamma} & \qin \Omega\,,\\[3mm]
            u=0 & \qon \partial \Omega\,.
        \end{array}
    \right .
    \label{hjhjss1}
\end{equation}

\begin{remark}\label{rem07}
Let us note that the definition of the function $\frac{f(x)}{u^\gamma}$ on the right hand side of \eqref{hjhjss1} needs to be precised. In all the paper, we will intend that the function $F(x,s)=\frac{f(x)}{u^\gamma}$, defined in $\Omega\times[0,+\infty[$ with values in $[0,+\infty]$, is $F(x,0):=0$ on the set $\{x\in\Omega:f(x)=0\}$.
\end{remark}

In this Section we will prove that for every $p>1$ the problem \eqref{hjhjss1} admits a weak solution in the following sense.

\begin {Definition}\label{mimi}
A function $u\in
W^{1,p}_0(\Omega)$ is a {\it weak solution} of problem \eqref{hjhjss1} if for every open set $\omega\subset\!\subset \Omega$ there exists $c_\omega>0$ such that
$u\ge c_\omega$ a.e. in $\omega$ and
\begin{equation*}\label{hjhjssqq}
    \int_\Omega |\nabla u|^{p-2}\nabla u \cdot \nabla v\, dx = \int_\Omega \frac{f}{u^\gamma}\, v \, dx\,
\end{equation*}
for every $v\in W^{1,p}_0(\Omega)$.
\end{Definition}

\begin{Theorem}\label{pproblem}
For any fixed $p>1$ and $0<\gamma\leq 1$, and for any $f\in L^N(\Omega)$, with $f\geq 0$, there exists a bounded unique weak solution of problem \eqref{hjhjss1} in the sense of Definition \ref{mimi}.
\end{Theorem}

\begin{pf}
By Theorems 4.1 and 4.4 of \cite{DCA} there exists $u\in
W^{1,p}_0(\Omega)\cap L^\infty(\Omega)$ distributional solution to
problem \eqref{hjhjss1}, i.e. a function $u$ such that $u\ge c_\omega$ a.e. in $\omega$ for every open set $\omega\subset\!\subset \Omega$ and
\begin{equation}\label{hjhjssff}
    \int_\Omega |\nabla u|^{p-2}\nabla u \cdot \nabla \varphi\, dx = \int_\Omega \frac{f}{u^\gamma}\, \varphi \,dx\,
\end{equation}
for every $\varphi\in C^{\infty}_0(\Omega)$. Given $v\in W^{1,p}_0(\Omega)$,
we will prove that
\begin{equation}\label{asda}
    \int_\Omega |\nabla u|^{p-2}\nabla u \cdot \nabla v\, dx = \int_\Omega \frac{f}{u^\gamma}\, v \, dx\,.
\end{equation}
%{{\color{red}
Observe that we may assume $v\ge0$ without loss of generality.

Consider a sequence $\varphi_n\in C^{\infty}_0(\Omega)$ such that $\varphi_n\ge0$ and
\begin{equation}\label{rrrr}
\varphi_n\to v\quad { \hbox { strongly in } } W^{1,p}_0(\Omega).
\end{equation}

%}}
We take for every $\eta>0$ the function
\begin{equation*}\label{hjhjssdd}
\rho_\eta*(v\wedge \varphi_n),
\end{equation*}
where $\rho_\eta$ is a standard convolution kernel and  $v\wedge \varphi_n:=\inf\{v,\varphi_n\}$. By taking it as test function in \eqref{hjhjssff}, we get
\begin{equation}\label{hjhjssffbb}
    \int_\Omega |\nabla u|^{p-2}\nabla u \cdot \nabla (\rho_\eta*(v\wedge \varphi_n))\, dx = \int_\Omega \frac{f}{u^\gamma}\, (\rho_\eta*(v\wedge \varphi_n))\,dx\,.
\end{equation}
We want to pass to the limit as $\eta\to 0$ and we get
\begin{equation*}\label{rrrrss}
\rho_\eta*(v\wedge \varphi_n)\to v\wedge \varphi_n\quad { \hbox { strongly in } } W^{1,p}_0(\Omega).
\end{equation*}
This implies for the left hand side of \eqref{hjhjssff} that
\begin{equation}\label{hjhjssffww}
    \int_\Omega |\nabla u|^{p-2}\nabla u \cdot \nabla (\rho_\eta*(v\wedge \varphi_n))\, dx
    \to
        \int_\Omega |\nabla u|^{p-2}\nabla u \cdot \nabla (v\wedge \varphi_n)\, dx\,.
    \end{equation}
    As far as the right hand side of \eqref{hjhjssffbb}, let us observe that for any $\eta>0$
$$\sop(\rho_\eta*(v\wedge \varphi_n))\subseteq K_n,$$
where $K_n$ is a compact set contained in $\Omega$, and that we have, choosing $\varphi\equiv1$ on $K_n$ in \eqref{hjhjssff},
\begin{equation*}
\frac{f}{u^\gamma}\in L^1(K_n).
\end{equation*}
Moreover, for any $\eta>0$
\begin{equation*}
\|\rho_\eta*(v\wedge \varphi_n)\|_{L^\infty}\leq
\|v\wedge \varphi_n\|_{L^\infty}
\end{equation*}
and, as $\eta\to 0$,
\begin{equation*}
\rho_\eta*(v\wedge \varphi_n)\to
v\wedge \varphi_n\qquad a.e.,
\end{equation*}
and so
\begin{equation*}
\rho_\eta*(v\wedge \varphi_n)\to
v\wedge \varphi_n\qquad w^*-L^\infty.
\end{equation*}
We conclude that, as $\eta\to 0$,
    \begin{equation}\label{floflo}
 \int_\Omega \frac{f}{u^\gamma}\, (\rho_\eta*(v\wedge \varphi_n))\,dx\to  \int_\Omega \frac{f}{u^\gamma}\, (v\wedge \varphi_n)\,dx.
\end{equation}
  By \eqref{hjhjssffbb},   \eqref{hjhjssffww} and \eqref{floflo}   we obtain
      \begin{equation}\label{flafla}
\int_\Omega |\nabla u|^{p-2}\nabla u \cdot \nabla (v\wedge \varphi_n)\, dx
=
 \int_\Omega \frac{f}{u^\gamma}\, (v\wedge \varphi_n)\,dx.
\end{equation}
    Now, we are going to pass to the limit in \eqref{flafla}, as $n\to +\infty$. Since
    \begin{equation*}
v\wedge \varphi_n\to
v\qquad \hbox{ in }   W^{1,p}(\Omega),
\end{equation*}
we have
          \begin{equation}\label{fluflu}
\int_\Omega |\nabla u|^{p-2}\nabla u \cdot \nabla (v\wedge \varphi_n)\, dx
\to
\int_\Omega |\nabla u|^{p-2}\nabla u \cdot \nabla v\, dx
.
\end{equation}
We are going now to prove that
\begin{equation*}\label{flifli}
\int_\Omega \frac{f}{u^\gamma}\, (v\wedge \varphi_n)\,dx
\to
\int_\Omega \frac{f}{u^\gamma}\, v\,dx\,,
\end{equation*}
as $n\to +\infty$. Indeed
\begin{equation*}\label{flicflic}
 \frac{f}{u^\gamma}\, (v\wedge \varphi_n)
\to
\frac{f}{u^\gamma}\, v\,\qquad \hbox{a.e.}
\end{equation*}
and
$$
0\leq \frac{f}{u^\gamma}\, (v\wedge \varphi_n)
\leq
\frac{f}{u^\gamma}\, v.
$$
Then by Lebesgue dominated convergence theorem, it is sufficient to prove that
\begin{equation}\label{flefle}
\frac{f}{u^\gamma}\,v \in L^1(\Omega)\,.
\end{equation}
We will prove that there exists a constant $C>0$, independent of $n$, satisfying
 for every $n\in \N$
\begin{equation}\label{estima12}
\int_\Omega \frac{f}{u^\gamma}|\varphi_n|\,dx\leq C\,.
\end{equation}
Indeed, by \eqref{hjhjssff} we have
\begin{multline*}
\int_\Omega \frac{f}{u^\gamma}\varphi_n\,dx=\int_\Omega |\nabla u|^{p-2}\nabla u\cdot \nabla \varphi_n \,dx\leq \int_\Omega |\nabla u|^{p-1} |\nabla \varphi_n| \,dx\\
\leq \int_\Omega |\nabla u|^p\,dx+\int_\Omega|\nabla \varphi_n|^{p} \,dx\le
C
\,,
\end{multline*}
where the last inequality is due to \eqref{rrrr} and to the
fact that $u\in
W^{1,p}_0(\Omega)$. This prove \eqref{estima12}.
By \eqref{rrrr}, \eqref{estima12}  and Fatou$^\prime$s Lemma, \eqref{flefle} follows.
Therefore formula \eqref{asda} holds.

Let us now prove the uniqueness of the solution. Assuming that there exist two solutions $u^1,u^2$, the uniqueness follows taking $u^1-u^2$ as test function in the two equations satisfied by $u^1$ and $u^2$ and observing that the term $\frac{f}{u^\gamma}$ is nonincreasing in the $u$ variable.
\end{pf}

\section{Main result}

This section is devoted to solve problem
\begin{equation}
    \left \{
        \begin{array}{cl}
            \displaystyle-\Div \left( \frac{Du}{|Du|}\right) = \frac{f(x)}{u^\gamma} & \qin \Omega\,,\\
            u=0 & \qon \partial \Omega\,,
        \end{array}
    \right .
    \label{problema1}
\end{equation}
for nonnegative data $f\in L^N(\Omega)$ and $0<\gamma\leq 1$. We begin by introducing the notion of solution to this problem.
%
%\begin{remark}\label{rem07}
%Let us note that the definition of the function $\frac{f(x)}{u^\gamma}$ on the right hand side of \eqref{problema1} needs to be precised. In all the paper, we will intend that the function $F(x,s)=\frac{f(x)}{u^\gamma}$, defined in $\Omega\times[0,+\infty[$ with values in $[0,+\infty]$, is $F(x,0):=0$ on the set $\{x\in\Omega:f(x)=0\}$.
%\end{remark}

\begin{Definition}\label{def1}
Let $f\in L^{N}(\Omega)$ with $f\ge0$ a.e..  We say that $u\in BV(\Omega)\cap L^\infty(\Omega)$, $u\geq 0$ a.e., is a \textbf{weak solution} of problem \eqref{problema1}
if there exists $\z\in\DM_{\rm{loc}}(\Omega)$ with $\|\z\|_{\infty} \le 1$ such that
\begin{enumerate}[$(a)$]
\item $\displaystyle \frac f{u^\gamma} \in L^1_{loc}(\Omega)\,,$

\item $\displaystyle\chi_{\{u>0\}}\in BV_{\rm{loc}}(\Omega)$\,,

\item  $\displaystyle-(\Div \z)\chi^*_{\{u>0\}}  =\frac f{u^\gamma} \,,\quad\text{ in }\, \dis (\Omega)\,,$

\item $\displaystyle(\z, Du)=|Du|\,, \ \  (\z, D\chi_{\{u>0\}})=|D\chi_{\{u>0\}}| \,\text{ as measures in }\, \Omega\,,$

\item    $u+ [u\z,\nu]=0 \ \mathcal H^{N-1}\text{--a.e. on }\, \partial\Omega\,.$

\end{enumerate}
\end{Definition}

\begin{remark}\label{rem0}\rm
Since
\[
\chi_{\{u>0\}}\in BV_{\rm{loc}}(\Omega),
\]
it follows from  Proposition \ref{nuova} that $\z\chi_{\{u>0\}}\in \DM_{\rm{loc}}(\Omega)$
and
the following equality holds
  \begin{equation}\label{divergence}
\Div\big(\z\chi_{\{u>0\}}\big)=\big(\Div \z\big)\chi^*_{\{u>0\}}+(\z,D \chi_{\{u>0\}}).
\end{equation}
Then, since $(\z, D\chi_{\{u>0\}})=|D\chi_{\{u>0\}}|$, the equation (c) is equivalent to
\begin{equation*}
\label{SSS}
-\Div\big(\z\chi_{\{u>0\}}\big)+|D \chi_{\{u>0\}}|=\frac f{u^\gamma}.
\end{equation*}
We recall that the term $|D \chi_{\{u>0\}}|$ is a measure concentrated on the reduced boundary $\partial^*\{u>0\}$ (or equivalently the reduced boundary $\partial^*\{u=0\}$). Moreover $|D \chi_{\{u>0\}}|$ coincides locally with the perimeter of the $\mathcal H^{N-1}$-rectifiable set $\{u=0\}$, i.e.
for every open set $\omega\subset\!\subset\Omega$ we have that $|D \chi_{\{u>0\}}|(\omega)$ coincides with the perimeter of $\{u=0\}\cap\omega$.
\end{remark}

\begin{remark}\label{rem1}\rm Let us point out that, by
$$
\int_\Omega\frac f{u^\gamma} |\varphi|\, dx<+\infty\,,\quad\forall \varphi\in C_0^\infty(\Omega)\,,
$$
any solution to \eqref{problema1} satisfies
$$
|\{x\in\Omega: u(x)=0, f(x)>0\}|=0\,,
$$
which means that
\begin{equation}\label{null1}
  \{x\in\Omega: u(x)=0\}\subseteq \{x\in\Omega: f(x)=0\}\,,
\end{equation}
except for a set of zero Lebesgue measure. Note that this implies that if $f>0$ a.e. in $\Omega$, then $u>0$ a.e. in $\Omega$. Moreover, if $f\not\equiv 0$, then for every solution $u$ we have $u\not\equiv 0$ and if $f\equiv 0$, then $u\equiv 0$ is a solution.

A further remark is in order. It follows from \eqref{null1} and by the fact that $(x,s)\mapsto \frac{f(x)}{s^\gamma}$ is a Carath\'eodory function on $\Omega\times [0,+\infty[$ with values in $[0,+\infty]$
(see Remark \ref{rem07}),
that
\begin{equation}\label{null2}
  \displaystyle \frac f{u^\gamma}=\frac f{u^\gamma}\chi_{\{u>0\}}\quad \hbox{a.e.}
\end{equation}
\end{remark}

\begin{remark}\label{rem2}\rm
We explicitly observe that there is a variational formulation of solution to \eqref{problema1}, see \eqref{var} below.
This formulation looks like that of renormalized solution.
\end{remark}

\begin{Theorem}\label{teoexist}
For every nonnegative $f\in L^{N}(\Omega)$ with $f\geq 0$ there is a weak solution to problem \eqref{problema1}.

Furthermore, for every nondecreasing and Lipschitz--continuous function $h:[0,+\infty[\to [0,+\infty[$ such that $h(0)=0$ and for every $\varphi\in C^\infty_0(\Omega)$,
\begin{equation}\label{var}
  \int_\Omega\varphi|Dh(u)|+\int_\Omega h(u)\z\cdot\nabla\varphi\, dx=\int_\Omega \frac f{u^\gamma}h(u)\varphi\, dx\,.
\end{equation}
Moreover  \begin{equation}\label{rrrrr}
h(u)+[h(u)\z,\nu]=0
\end{equation}
holds $\mathcal H^{N-1}$--a.e. on $\partial\Omega $.
\end{Theorem}

\begin{pf}
The proof will be divided in several steps.

{\sl Step 1. Approximating problems.}

For every $p>1$ let us consider the following problem

\begin{equation}
    \left \{
        \begin{array}{cl}
            \displaystyle-\Delta_p \left(u_p\right) = \frac{f(x)}{u_p^\gamma} & \qin \Omega\,,\\
            u_p=0 & \qon \partial \Omega\,.
        \end{array}
    \right .
    \label{problemap}
\end{equation}

By Theorem \ref{pproblem} there exists $u_p\in
W^{1,p}_0(\Omega)\cap L^\infty(\Omega)$ weak solution to
problem \eqref{problemap}, i.e. it satisfies
\begin{equation}\label{hjhjss}
    \int_\Omega |\nabla u_p|^{p-2}\nabla u_p \cdot \nabla v\, dx = \int_\Omega \frac{f}{u_p^\gamma}\, v\, dx\,
\end{equation}
for every $v\in W^{1,p}_0(\Omega)$. Moreover for every open set $\omega\subset\!\subset \Omega$ there exists $c_\omega>0$ such that
$u_p\ge c_\omega$ a.e. in $\omega$.

  {\sl Step 2. $BV$--estimate.}

Taking the test function $u_p$ in problem
\eqref{hjhjss},  and by using the H\"older inequality we get
\begin{equation}\label{BV}
\int_\Omega |\nabla u_p|^p dx= \int_\Omega fu^{1-\gamma}_p \, dx \le \|f\|_N |\Omega|^{\gamma\frac{N-1}{N}}\|u_p\|^{1-\gamma}_{\frac N{N-1}}.
\end{equation}
Thus, applying the Sobolev and the Young inequalities we have
\begin{equation}\label{BV1}
\|u_p\|_{\frac N{N-1}}\leq S \int_\Omega |\nabla u_p| dx\leq \frac Sp\int_\Omega |\nabla u_p|^p dx+
\frac{S(p-1)}{p}|\Omega|
\end{equation}
\[
\leq S\|f\|_N|\Omega|^{\gamma\frac{N-1}{N}}\|u_p\|_{\frac N{N-1}}^{1-\gamma}+S|\Omega|.
\]
Since $1-\gamma<1$, it follows that $(u_p)$ is bounded in $L^{\frac N{N-1}}(\Omega)$.
Therefore, by \eqref{BV} and \eqref{BV1},  we have
\begin{equation}\label{BV12}
\int_\Omega |\nabla u_p|^p dx\leq M,
\end{equation}
for certain constant which does not depend on $p$.
Thanks to Young's inequality, it implies that
\begin{equation*}\label{BV13}
\int_\Omega |\nabla u_p| dx\leq \frac1p\int_\Omega |\nabla u_p|^p dx+\frac{p-1}p|\Omega|\le M+|\Omega|\,.
\end{equation*}
Having in mind that $u_p\big|_{\partial\Omega}=0$, we obtain that
$(u_p)$ is bounded in $BV(\Omega)$. Hence,
there exists a function $u$ such that,
\begin{equation}\label{uBV}
u\in BV(\Omega),
\end{equation}
and, up to a
subsequence,
\begin{equation}
u_p \to u\ {\text{ in }}\  L^q(\Omega) \ {\text { for every }}\  1\leq q< {\frac N{N-1}}\,,\label{23BIS}
\end{equation}
\begin{equation}
u_p \rightharpoonup u  \ {\text {  \ \ in\ \  }}L^{\frac N{N-1}}(\Omega)\  \label{24BISBIS}
\end{equation}
and
$u_p \to u$   pointwise a.e. in $\Omega$\,.
Moreover
\begin{equation}
u\geq 0 \ {\text { a.e.\,, }}\label{convergence}
\end{equation}
and
$
\nabla u_p \weakly Du$ weakly as measures.
%
%up to a
%subsequence, $u_p \to u$ in $L^1(\Omega)$ (and so a.e. in $\Omega$) and $u_p \to u$ in $L^q(\Omega)$, with $1\leq q< {\frac N{N-1}}$.  Note that $u\geq 0$ a.e. in $\Omega$.
%and $D u_n$ converges to
%$Du$ $*$--weakly as measures when $n \to \infty$.

%
%Taking the function test $\frac{T_k(u_n)}{k}$ in problem
%\eqref{hjhjss}, we get
%\[
%    \frac{1}{k}\int_\Omega (\z_n,DT_k(u_n)) + \frac{1}{k}\int_\Omega T_k(u_n)^*|Du_n| = \int_\Omega f_n \frac{T_k(u_n)}{k}\, dx \le \int_\Omega f_n\, dx,,
%\]
%where $C$ does not depend on $n$. Since $(\z_n,Du_n) =|Du_n|$, it
%follows from Proposition \ref{prop1} that $(\z_n,DT_k(u_n))=|DT_k(u_n)|$, which is nonnegative. Thus
%\[
%    \frac{1}{k}\int_\Omega T_k(u_n)^*|Du_n| \le C\,.
%\]
%Then, letting $k \to \infty$ in the inequality above we arrive at
%\[
%    \int_\Omega |Du_n| \le C\,.
%\]
%Therefore, $u_n$ is bounded in $BV(\Omega)$ and, up to a
%subsequence, $u_n \to u$ in $L^1(\Omega)$ and $D u_n$ converges to
%$Du$ $*$--weakly as measures when $n \to \infty$.

  {\sl Step 3. $L^\infty$--estimate.}

Taking the test function $G_k(u_p)=(u_p-k)^+$, with $k\in ]0,+\infty[$\,, in problem
\eqref{hjhjss},  and by using the H\"older inequality we get
\begin{multline*}\label{Linfty}
\int_\Omega |\nabla G_k(u_p)|^p\, dx= \int_\Omega \frac{f}{u^\gamma_p}G_k(u_p) \, dx\\
 \le \frac1{k^\gamma}\int_\Omega fG_k(u_p) dx\leq
\frac{\|f\|_N}{k^\gamma}\|G_k(u_p)\|_{\frac N{N-1}}.
\end{multline*}

Thus, the Sobolev and the Young inequalities imply
\begin{multline*}
\|G_k(u_p)\|_{\frac N{N-1}}\leq S \int_\Omega |\nabla G_k(u_p)|\, dx\leq \frac Sp\int_\Omega |\nabla G_k(u_p)|^p\, dx+
\frac{S(p-1)}{p}|\Omega|\\
\leq S\frac{\|f\|_N}{k^\gamma}\|G_k(u_p)\|_{\frac N{N-1}}+S(p-1)|\Omega|.
\end{multline*}
Now, by choosing $k$ satisfying $S\frac{\|f\|_N}{k^\gamma}<1$,
we get
\[
\|G_k(u_p)\|_{\frac N{N-1}}\leq
C(p-1),
\]
with $C$ independent of $p$. By applying Fatou's Lemma we deduce
\[
\|G_k(u)\|_{\frac N{N-1}}\leq \liminf_{p\to 1}\|G_k(u_p)\|_{\frac N{N-1}}\leq  \liminf_{p\to 1}C(p-1)=0.
\]
We conclude that $0\leq u\leq k$ a.e. in $\Omega$ and so
\begin{equation*}\label{uLinfty}
u\in L^\infty(\Omega).
\end{equation*}

  {\sl Step 4. $L^1_{\rm loc}$--estimate for the singular term.}

We are going to prove in this step that there exists a constant $C>0$, independent of $p$, satisfying
\begin{equation}\label{estima}
\int_\Omega \frac{f}{u_p^\gamma}|\varphi|\,dx\leq C\,,\quad \forall \varphi\in C_0^\infty(\Omega),
\ \forall 1<p\leq p_0\,.
\end{equation}
Fixed $\varphi\in C_0^\infty(\Omega)$, with $\varphi\geq 0$, by \eqref{BV12} we have
\begin{multline*}
\int_\Omega \frac{f}{u_p^\gamma}\varphi\,dx=\int_\Omega |\nabla u_p|^{p-2}\nabla u_p\cdot \nabla \varphi \,dx\leq \int_\Omega |\nabla u_p|^{p-1} |\nabla \varphi| \,dx\\
\leq \int_\Omega |\nabla u_p|^p\,dx+\int_\Omega|\nabla \varphi|^{p} \,dx\le
M+\int_\Omega\Big(|\nabla \varphi|^{p_0}+1\Big) \,dx<+\infty\,.
\end{multline*}
This implies that \eqref{estima} holds for every $\varphi\in C_0^\infty(\Omega)$, even if it changes sign.
By this estimate and Fatou's Lemma we get
\begin{equation}\label{estima1}
\int_\Omega \frac{f}{u^\gamma}|\varphi|\,dx<+\infty\,,\ \ \ \forall \varphi\in C_0^\infty(\Omega)\,.
\end{equation}

  {\sl Step 5. Vector field $\z$.}

Now, we want to find a vector field $\z\in \DM_{\rm{loc}}(\Omega)$ with $\|\z\|_{\infty} \le 1$, to play the role of $\frac{Du}{|Du|}$.
In this Step, we will follow the argument of \cite{MRST}.

Fix $1\le q<\infty$ and consider $1<p<q^\prime$. It follows from \eqref{BV12}, that
\begin{equation*}
\int_\Omega \Big||\nabla u_p|^{p-2}\nabla u_p\Big|^q dx= \int_\Omega |\nabla u_p|^{q(p-1)}  \,dx
\end{equation*}
\[
\leq \left(\int_\Omega |\nabla u_p|^p\,dx\right)^{\frac{q}{p^\prime}}|\Omega|^{1-\frac q{p^\prime}}\leq M^{\frac{q}{p^\prime}}|\Omega|^{1-\frac q{p^\prime}}.
\]
Hence,
\begin{equation*}\label{estima2}
\||\nabla u_p|^{p-2}\nabla u_p\|_q\le M^{\frac{1}{p^\prime}}|\Omega|^{\frac1q-\frac 1{p^\prime}}\le(1+M)^{\frac1q}(1+|\Omega|)^{\frac1q}
\end{equation*}
and so the sequence $(|\nabla u_p|^{p-2}\nabla u_p)_p$ is bounded in $L^q(\Omega;\R^N)$. Then there exists $\z_q\in L^q(\Omega;\R^N)$ such that, up to subsequences,
\[
|\nabla u_p|^{p-2}\nabla u_p\rightharpoonup \z_q\,,\quad\hbox{weakly in }L^q(\Omega;\R^N)\,.
\]
A diagonal argument shows that there exists a vector field $\z$ (independent of $q$) satisfying
\begin{equation}\label{convergence2}
|\nabla u_p|^{p-2}\nabla u_p\rightharpoonup \z\,,\quad\hbox{weakly in }L^q(\Omega;\R^N)\,,\quad\forall 1\leq q<\infty\,.
\end{equation}
Moreover, by applying the lower semicontinuity of the $q$--norm, the previous estimate $\||\nabla u_p|^{p-2}\nabla u_p\|_q\le M^{\frac{1}{p^\prime}}|\Omega|^{\frac1q-\frac 1{p^\prime}}$ implies
\[
\|\z\|_q\le |\Omega|^{\frac1q}\,,\quad\forall q<\infty\,,
\]
so that, letting $q$ go to $\infty$, we obtain $\z\in L^\infty(\Omega;\R^N)$ and $\|\z\|_\infty\le1$.

Using $\varphi \in C^\infty_0 (\Omega)$, with $\varphi \ge 0$, as a function test in \eqref{hjhjss}, we arrive at
\begin{equation}\label{hjhj}
    \int_\Omega |\nabla u_p|^{p-2}\nabla u_p \cdot \nabla \varphi\, dx = \int_\Omega \frac{f}{u_p^\gamma}\, \varphi\, dx\,,
\end{equation}
and when we take $p\to 1$, using \eqref{convergence}, Fatou's Lemma and \eqref{convergence2} it becomes
\[
\int_\Omega \z \cdot \nabla\varphi\, dx  \geq \int_\Omega \frac{f}{u^\gamma} \,\varphi\, dx\,.
\]
Therefore,
\begin{equation}\label{iigg}
    -\Div \z  \geq  \frac{f}{u^\gamma}\qin \dis (\Omega)\,,
\end{equation}
so $-\Div \z$ is a nonnegative Radon measure. By \eqref{estima} and \eqref{hjhj}
we have that
\[
0\le -\int_\Omega \varphi\, \Div\z=\int_\Omega \z\cdot\nabla\varphi\, dx  <+\infty
\]
holds for every $\varphi \in C^\infty_0 (\Omega)$ satisfying $\varphi \ge 0$.
This implies that the total variation of $-\Div \z$ is locally finite. Therefore
$\z\in \DM_{\rm {loc}}(\Omega)$.

Let us note that the total variation of $-\Div \z$ is only locally bounded, since by \eqref{estima} the term $\frac{f}{u_p^\gamma}$ is bounded only in $L^1_{\rm {loc}}(\Omega)$ (see also \eqref{hjhj}).

{\sl Step 6. The following equation
\begin{equation}
\label{equ:main}  -u^*\, \Div\z= f u^{1-\gamma}
\end{equation}
 holds in $\dis (\Omega)$ and the function $u^*$ belongs to $L^1(\Omega,\Div\z)$.}

Firstly, let us prove that by \eqref{iigg} we have
\begin{equation}\label{mlml}
-u^*\, \Div\z\geq f u^{1-\gamma}\,,\quad \hbox{in }\dis (\Omega),
  \end{equation}
namely
\begin{equation}
  -\int_\Omega u^*\varphi\, \Div\z\geq \int_\Omega f u^{1-\gamma}\varphi\,\quad\forall \varphi\in C^\infty_0 (\Omega).
  \end{equation}

%{\color{red}{
Indeed we can apply \eqref{iigg} using the test function $(u*\rho_\epsilon)\varphi$, where $\rho_\epsilon$ is a standard mollifier and $\varphi\in C^\infty_0 (\Omega)$.
We get
\begin{equation}\label{mlml1}
 -\int_\Omega (u*\rho_\epsilon)\varphi\, \Div\z\geq \int_\Omega  \frac{f}{u^\gamma} \,(u*\rho_\epsilon)\varphi\,dx.
  \end{equation}
We recall that, since $u\in L^\infty(\Omega)$ then we have also that $|u^*|\leq \|u\|_\infty$ $\mathcal H^{N-1}$-a.e. This implies that $|u^*|\leq \|u\|_\infty$ $\Div \z$-a.e. since $\Div\z<\!\!<\mathcal H^{N-1}$ (this gives also that $u^*\in L^1_{\rm{loc}}(\Omega , \Div \z)$\,). By Proposition 3.64 (b) and Proposition 3.69 (b) of \cite{AFP} we have that
  \begin{equation*}
  u*\rho_\epsilon \to u^*   \quad\mathcal H^{N-1}-{\hbox {a.e.}}
    \end{equation*}
   and therefore $\Div \z$-a.e.
    We can pass to the limit in both sides of the inequality \eqref{mlml1} by dominated convergence theorem (with respect to the measure $\Div\z$ on the left hand side and the Lebesgue measure on the right hand side), since
$|u*\rho_\epsilon|\leq \|u\|_\infty$ in $\Omega$ and since $\frac{f}{u^\gamma}\varphi\in L^1(\Omega)$ by \eqref{estima1}. Therefore \eqref{mlml} holds.
%}}

Now, in order to prove the opposite inequality, let $\varphi \in C^\infty_0 (\Omega)$, with $\varphi \ge 0$, and take $u_p\varphi$ as a function test in \eqref{hjhjss}. Then
\[
\int_\Omega\varphi|\nabla u_p|^p\, dx+\int_\Omega u_p|\nabla u_p|^{p-2}\nabla u_p\cdot\nabla\varphi\, dx=\int_\Omega f u_p^{1-\gamma}\varphi\, dx\,.
\]
Applying Young's inequality, it yields
\begin{multline*}
  \int_\Omega\varphi|\nabla u_p|\, dx+\int_\Omega u_p\,|\nabla u_p|^{p-2}\nabla u_p\cdot\nabla\varphi\, dx \\
  \le \frac1p \int_\Omega\varphi|\nabla u_p|^p\, dx+\int_\Omega u_p\,|\nabla u_p|^{p-2}\nabla u_p\cdot\nabla\varphi\, dx+\frac{p-1}p\int_\Omega \varphi\, dx\\
  \le \int_\Omega f u_p^{1-\gamma}\varphi\, dx+\frac{p-1}p\int_\Omega \varphi\, dx\,.
\end{multline*}
%{{\color{red}
To pass to the limit on the left hand side, we apply the lower semicontinuity, jointly with \eqref{23BIS} and \eqref{convergence2}. Assume first that $0<\gamma<1$.
On the right hand side, we point out that, as $p\to 1$, by \eqref{24BISBIS}
\begin{equation}\label{GAMMA}
u_p^{1-\gamma}\rightharpoonup u^{1-\gamma}\ \  \hbox{ in \ } L^{\frac{N}{N-1}\frac{1}{1-\gamma}}(\Omega)\subseteq L^{\frac{N}{N-1}}(\Omega)\,,
\end{equation}
so that
\begin{equation}\label{stella}
\lim_{p\to 1}\int_\Omega f u_p^{1-\gamma}\varphi\, dx=\int_\Omega f u^{1-\gamma}\varphi\, dx\,.
\end{equation}
It leads to
\begin{equation}\label{inequ:main0}
\int_\Omega\varphi|Du|+\int_\Omega u\z\cdot\nabla\varphi\, dx
\le \int_\Omega f u^{1-\gamma}\varphi\, dx\,.
\end{equation}
Note that if $\gamma=1$ the first integral on the right hand side does note depend on $p$.
Taking into account \eqref{ALPHA} and \eqref{ec:2},
we deduce by \eqref{inequ:main0} that
\begin{multline*}
  -\int_\Omega u^*\varphi\, \Div\z= \int_\Omega\varphi (\z, Du)+\int_\Omega u\z\cdot\nabla \varphi\, dx\\
  \le \int_\Omega\varphi \, |Du|+\int_\Omega u\z\cdot\nabla \varphi\, dx
  \le \int_\Omega f u^{1-\gamma}\varphi\, dx\,.
\end{multline*}
%}}
Therefore,
\begin{equation}\label{BETA}
-u^*\, \Div\z\le f u^{1-\gamma}\,,\quad \hbox{in }\dis (\Omega)\,
\end{equation}
%wherewith
%\[
%-\big( \Div\z\big)\chi_{\{u^*>0\}}\le \frac{f}{u^{\gamma}}\,,\quad \hbox{in }\dis (\Omega)\,,
%\]
%So, it follows from here and \eqref{dis1} that $-\big(\Div \z\big)\chi_{\{u^*>0\}}  =  \frac{f}{u^\gamma}$ holds in $\dis (\Omega)$.
%It is worth pointing out that we have also shown
and by \eqref{mlml} and \eqref{BETA} we conclude that \eqref{equ:main} holds.
This implies that $  -u^*\, \Div\z\in L^1(\Omega)$, so that $u^*\in L^1(\Omega , \Div\z)$.

  {\sl Step 7. The equality $(\z,Du)=|Du|$, as measures on $\Omega$, holds.}

It is a straightforward consequence of \eqref{inequ:main0} and \eqref{equ:main}. Indeed, consider
 $\varphi \in C^\infty_0 (\Omega)$ with $\varphi \ge 0$. Then we have
\begin{equation*}\label{BVBBVV}
\int_\Omega \varphi |D u| +\int_\Omega u\, \z\cdot\nabla \varphi\, dx\leq
\int_\Omega fu^{1-\gamma}\varphi \, dx=-\int_\Omega u^* \varphi\, \Div\z \,.
\end{equation*}
Therefore by Proposition \ref{nuova}
\begin{equation*}\label{BVBBVVV}
\int_\Omega\varphi |D u| \leq-\int_\Omega u\, \z\cdot\nabla \varphi\, dx-
\int_\Omega u^*\varphi \,\Div\z = \int_\Omega \varphi(\z,Du).
\end{equation*}
The arbitrariness of $\varphi$ implies that
\begin{equation*}\label{BVBBVVVB}
|D u| \leq (\z,Du)\,,
\end{equation*}
as measures on $\Omega$.
On the other hand, since $\|z\|_\infty\leq 1$, the opposite inequality holds, i.e.
\begin{equation*}\label{BVBBVVVB}
|D u| \geq (\z,Du)\,,
\end{equation*}
as measures on $\Omega$.
This concludes the proof of Step 7.

  {\sl Step 8. The boundary condition $   u+ [u\z,\nu]=0$ holds $\mathcal H^{N-1}\text{--a.e. on }\, \partial\Omega $.}

 By taking $u_p$ as test function in \eqref{hjhjss} we have
\begin{equation*}\label{BVBBM}
\int_\Omega |\nabla u_p|^p  \,dx=
\int_\Omega fu^{1-\gamma}_p\, dx \,.
\end{equation*}
On the other hand, by the fact that $u_p=0$ on $\partial \Omega$ and by the Young inequality we get
\begin{equation*}\label{BVBBVS}
\begin{split}
&\int_\Omega |\nabla u_p| \,dx+\int_{\partial\Omega} u_p \,d\mathcal H^{N-1}
\leq\frac1p\int_\Omega |\nabla u_p|^p \,dx+\frac{p-1}{p}|\Omega| \\
\leq&\int_\Omega fu^{1-\gamma}_p \, dx+\frac{p-1}{p}|\Omega| \,.
\end{split}
\end{equation*}
 We may use the lower semicontinuity of the functional on the left hand side to pass to the limit as $p\to 1$ and obtain by \eqref{GAMMA} and \eqref{equ:main}
\begin{equation*}\label{BVBBVVC}
\int_\Omega |D u| +\int_{\partial\Omega} u\, d\mathcal H^{N-1}\leq
\int_\Omega fu^{1-\gamma} \, dx= -\int_{\Omega} u^*\Div\z\,.
\end{equation*}
Applying the Green formula (see Proposition 2.6)
on the right hand side we have
\begin{equation*}\label{BVBBVVCXX}
\int_\Omega |D u| +\int_{\partial\Omega} u\, d\mathcal H^{N-1}\leq
\int_\Omega (\z,Du)-\int_{\partial\Omega} [u\z,\nu]\, d\mathcal H^{N-1}
\,.
\end{equation*}
Then, by Step 7, we arrive at
\[
\int_{\partial\Omega} (u+[u\z,\nu])\, d\mathcal H^{N-1}\leq0\,.
\]
Since $|[u\z, \nu]|\le u\|\z\|_\infty\le u$, we conclude that
\begin{equation*}\label{BVBBVVCXXM}
u+[u\z,\nu] =
0\quad \mathcal H^{N-1}\hbox{--a.e. on }\partial\Omega\,,
\end{equation*}
which concludes the proof of the Step 8.

   {\sl Step 9. Variational formulation.}

 Consider a Lipschitz--continuous and nondecreasing function $h:[0,+\infty[\to [0,+\infty[$ such that $h(0)=0$. Then by \eqref{uBV} we have that $h(u)\in BV(\Omega)$.
 Following the argument of Step 6, we have
 \begin{equation}\label{www}
-h(u)^*\, \Div\z\geq  \frac{f}{u^{\gamma}}h(u)\,,\quad \hbox{in }\dis (\Omega).
 \end{equation}

 In order to prove the opposite inequality,
let $\varphi \in C^\infty_0 (\Omega)$, with $\varphi \ge 0$, and taking $h(u_p)\varphi$ as test function  in \eqref{hjhjss} we obtain

%%%%%%%%
\[
\int_\Omega\varphi|\nabla u_p|^ph^\prime(u_p)\, dx+\int_\Omega h(u_p)|\nabla u_p|^{p-2}\nabla u_p\cdot\nabla\varphi\, dx=\int_\Omega \frac f{u^\gamma}h(u_p)\varphi\, dx\,.
\]
Applying Young's inequality, it yields
\begin{multline}\label{RR}
  \int_\Omega\varphi|\nabla h(u_p)|\, dx+\int_\Omega h(u_p)\,|\nabla u_p|^{p-2}\nabla u_p\cdot\nabla\varphi\, dx
\\
=  \int_\Omega\varphi|\nabla u_p|h^\prime(u_p)\, dx+\int_\Omega h(u_p)\,|\nabla u_p|^{p-2}\nabla u_p\cdot\nabla\varphi\, dx \\
  \le \frac1p \int_\Omega\varphi|\nabla u_p|^ph^\prime(u_p)\, dx+\int_\Omega h(u_p)\,|\nabla u_p|^{p-2}\nabla u_p\cdot\nabla\varphi\, dx+\frac{p-1}p\int_\Omega h^\prime(u_p)\varphi\, dx\\
  \le \int_\Omega \frac f{u_p^\gamma}h(u_p)\varphi\, dx+\frac{p-1}pL \int_\Omega \varphi\, dx\,,
\end{multline}
where $L$ denotes the Lipschitz constant of $h$.
Moreover, it follows from
\begin{equation}\label{stellastella}
h(u_p)\leq L u_p\,,
\end{equation}
the convergence \eqref{23BIS} and Vitali's Theorem that, as $p\to 1$,
\begin{equation}\label{SS}
h(u_p)\to h(u),\ {\text{ in }}\  L^q(\Omega) \ {\text { for every }}\  1\leq q< {\frac N{N-1}}\,.
\end{equation}
By \eqref{stellastella} we also get
\begin{equation*}\label{stellastella1}
\int_E \frac f{u_p^\gamma}h(u_p)\varphi\, dx\leq L
\int_E  f{u_p^{1-\gamma}}\varphi\, dx
\end{equation*}
for every Borel measurable set $E\subset\Omega$. Therefore, by \eqref{stella} and Vitali$^\prime$s Theorem we obtain
\begin{equation}\label{stellastellastella}
\lim_{p\to 1}\int_\Omega \frac f{u_p^\gamma}h(u_p)\varphi\, dx=\int_\Omega \frac f{u^\gamma}h(u)\varphi\, dx\,.
\end{equation}
The lower semicontinuity, jointly with \eqref{RR}, \eqref{SS}  and \eqref{stellastellastella}, leads to
\begin{equation}\label{inequ:main}
\int_\Omega\varphi|Dh(u)|+\int_\Omega h(u)\z\cdot\nabla\varphi\, dx
\le \int_\Omega \frac f{u^\gamma}h(u)\varphi\, dx\,.
\end{equation}
By \eqref{ALPHA}
we deduce that
\begin{multline*}
  -\int_\Omega h(u)^*\, \Div\z\,\varphi= \int_\Omega\varphi (\z, Dh(u))+\int_\Omega h(u)\z\cdot\nabla \varphi\, dx\\
  \le \int_\Omega\varphi|Dh(u)|+\int_\Omega h(u)\z\cdot\nabla\varphi\, dx
\le \int_\Omega \frac f{u^\gamma}h(u)\varphi\, dx
\end{multline*}
which implies that
\begin{equation}\label{inequ:main33}
-h(u)^*\, \Div\z\leq  \frac{f}{u^{\gamma}}h(u)\,,\quad \hbox{in }\dis (\Omega).
\end{equation}
By \eqref{www} and \eqref{inequ:main33}, we  conclude that
\begin{equation}\label{equ:main1}
  -h(u)^*\, \Div\z= \frac f{u^\gamma}h(u)\,,\quad \hbox{in }\dis (\Omega)\,.
\end{equation}

Next, we deduce that
\begin{equation}\label{nonna}
(\z, Dh(u))=|Dh(u)|.
\end{equation}
We point out that we cannot apply \cite[Proposition 2.8]{An} since $u$ is not continuous (it is proved for pairings such as $\z\in\DM(\Omega)$ and $u\in BV(\Omega)\cap C(\Omega)\cap L^\infty(\Omega)$). Thus, we must focus on our equation and repeat the same arguments as in Step 7.
In fact, let us consider
 $\varphi \in C^\infty_0 (\Omega)$ with $\varphi \ge 0$. Then by \eqref{inequ:main} and \eqref{equ:main1} we have
\begin{equation*}\label{BVBBVV}
\int_\Omega \varphi |D h(u)| +\int_\Omega h(u)\, \z\cdot\nabla \varphi\, dx\leq
\int_\Omega \frac{f}{u^{\gamma}}h(u)\varphi \, dx=-\int_\Omega \varphi h(u)^*\Div\z \,.
\end{equation*}
Therefore by \eqref{dist1}
\begin{equation*}\label{BVBBVVV}
\int_\Omega\varphi |D h(u)| \leq-\int_\Omega h(u)\, \z\cdot\nabla \varphi\, dx-
\int_\Omega h(u)^*\varphi \,\Div\z = \int_\Omega \varphi(\z,Dh(u)).
\end{equation*}
The arbitrariness of $\varphi$ implies that
\begin{equation*}
|D h(u)| \leq (\z,Dh(u))\,,
\end{equation*}
as measures on $\Omega$.
On the other hand, since $\|\z\|_\infty\leq 1$, the opposite inequality holds, i.e.
\begin{equation*}
|D h(u)| \geq (\z,Dh(u))\,,
\end{equation*}
as measures on $\Omega$.
This proves that $(\z, Dh(u))=|Dh(u)|$.

By \eqref{nonna}, \eqref{equ:main1} and \eqref{ALPHA} we have
 \begin{equation}\label{renormali}
|D h(u)|-\Div\big(h(u)\z\big)= \frac f{u^\gamma}h(u)\,,\quad \hbox{in }\dis (\Omega)
\end{equation}
and then \eqref{var} holds.
%{{\color{red}
Next,  take $h(u_p)$ as test function  in \eqref{hjhjss} to obtain
\[
\int_\Omega|\nabla u_p|^ph^\prime(u_p)\, dx=\int_\Omega \frac f{u_p^\gamma}h(u_p)\, dx\,.
\]
The boundary condition and Young's inequality yield
\begin{multline}\label{ec:5}
  \int_\Omega|\nabla h(u_p)|\, dx+\int_{\partial\Omega} h(u_p)\, d\mathcal H^{N-1}
=  \int_\Omega|\nabla u_p|h^\prime(u_p)\, dx \\
  \le \frac1p \int_\Omega|\nabla u_p|^ph^\prime(u_p)\, dx+\frac{p-1}p\int_\Omega h^\prime(u_p)\, dx\\
  \le \int_\Omega \frac f{u_p^\gamma}h(u_p)\, dx+\frac{p-1}p \int_\Omega h^\prime(u_p)\, dx\,.
\end{multline}
Having in mind \eqref{stellastella}, we deduce
\[
\lim_{p\to1}\int_\Omega \frac f{u_p^\gamma}h(u_p)\, dx=\int_\Omega \frac f{u^\gamma}h(u)\, dx\,.
\]
On the left hand side of \eqref{ec:5}, one can use the lower semicontinuity of the functional to arrive at
\begin{equation*}
  \int_\Omega|D h(u)|+\int_{\partial\Omega} h(u)\, d\mathcal H^{N-1}\le \int_\Omega \frac f{u^\gamma}h(u)\, dx\,.
\end{equation*}
Applying \eqref{renormali} and Green's formula \eqref{GreenIII}, this inequality becomes
\begin{equation*}
  \int_\Omega|D h(u)|+\int_{\partial\Omega} h(u)\, d\mathcal H^{N-1}\le \int_\Omega |D h(u)|-\int_{\partial\Omega}[h(u)\z,\nu]\, d\mathcal H^{N-1}\,,
\end{equation*}
wherewith
\begin{equation*}
  \int_{\partial\Omega} h(u)+[h(u)\z,\nu]\, d\mathcal H^{N-1}\le0\,.
\end{equation*}
Hence, $h(u)+[h(u)\z,\nu]=0$ holds $\mathcal H^{N-1}$--a.e. on $\partial\Omega $, i.e. \eqref{rrrrr} is proved.

   {\sl Step  10. The function $\chi_{\{u>0\}}$ belongs to $BV_{\rm{loc}}(\Omega)$.}

 Let $(h_n)$ be a sequence of Lipschitz--continuous and nondecreasing functions $h_n:[0,+\infty[\to[0,+\infty[$ such that $h_n(0)=0$
 and $h_n(u)\uparrow \chi_{\{u>0\}}$  in $L^1(\Omega)$.
% {{\color{red}
 (For instance, we may consider the truncation $h_n(s)=nT_{1/n}(s^+)$.)
% }}

 By \eqref{renormali}, for every $n$  and $\varphi \in C^1_0 (\Omega)$, with $\varphi \ge 0$, we have
  \begin{equation*}\label{inequ:main122}
\int_\Omega\varphi|D h_n(u)|\leq -\int_\Omega h_n(u)\z\cdot\nabla\varphi\, dx
+\int_\Omega \frac f{u^\gamma}h_n(u)\varphi\, dx\,.
\end{equation*}
On account of \eqref{estima1}, this fact implies that  for every $\omega\subset\subset\Omega$ there exists a constant $C$ such that
  \begin{equation*}\label{inequ:main1224}
|D h_n(u)|(\omega)\leq C\,.
\end{equation*}
Hence for every $\omega\subset\subset\Omega$
  \begin{equation*}\label{inequ:main1223}
|D \chi_{\{u>0\}}|(\omega)\leq\liminf_{n\to+\infty}|D h_n(u)|(\omega)\leq C\,,
\end{equation*}
and $\chi_{\{u>0\}}$ belongs to $BV_{\rm{loc}}(\Omega)$. Then by Proposition \ref{nuova} we have that $\z\chi_{\{u>0\}}\in \DM_{\rm{loc}}(\Omega)$
and by \eqref{ALPHA}
the following equality holds

  \begin{equation*}\label{divergence}
\Div\big(\z\chi_{\{u>0\}}\big)=\big(\Div \z\big)\chi^*_{\{u>0\}}+(\z,D \chi_{\{u>0\}})
\end{equation*}
 in $\dis (\Omega).$

   {\sl Step  11. $u$ satisfies the inequality $-\Div\big(\z\chi_{\{u>0\}}\big)+|D \chi_{\{u>0\}}|\leq\frac f{u^\gamma}\,,\quad \hbox{in }\dis (\Omega)\,.$}

By \eqref{renormali} for every $n$  and $\varphi \in C^1_0 (\Omega)$, with $\varphi \ge 0$, we have
 \begin{equation*}
\int_\Omega \varphi|D h_n(u)|+\int_\Omega h_n(u)\z\cdot\nabla\varphi\,dx=\int_\Omega \frac f{u^\gamma}h_n(u)\varphi\,dx\,.
\end{equation*}
By the lower semicontinuity we have
 \begin{equation*}
\int_\Omega \varphi|D \chi_{\{u>0\}}|+\int_\Omega \chi_{\{u>0\}}\z\cdot\nabla\varphi\,dx\leq\int_\Omega \frac f{u^\gamma}\varphi\,dx\,.
\end{equation*}

   {\sl Step  12. $u$ satisfies the equation $-\big(\Div \z\big)\chi^*_{\{u>0\}}  =  \frac{f}{u^\gamma}$ in $\dis (\Omega)$.}

Since
$
\chi_{\{u>0\}}\in BV_{\rm{loc}}(\Omega),
$ we have that \eqref{null2}  and \eqref{iigg} imply
\begin{equation*}\label{dis1}
    -\big(\Div \z\big)\chi^*_{\{u>0\}}  \geq  \frac{f}{u^\gamma}\qin \dis (\Omega)\,.
  \end{equation*}
Since $\|z\|_\infty\leq 1$, the inequality
\begin{equation}\label{uuu}
(\z, D\chi_{\{u>0\}})\leq|D\chi_{\{u>0\}}|
\end{equation}
holds. Moreover, by \eqref{divergence}, \eqref{uuu} and Step 11  we have that
\begin{equation*}
\begin{split}
-\big(\Div \z\big)\chi^*_{\{u>0\}}  &=-\Div\big(\z\chi_{\{u>0\}}\big)+(\z,D \chi_{\{u>0\}})\\
&\leq-\Div\big(\z\chi_{\{u>0\}}\big)+|D \chi_{\{u>0\}}|\leq\frac f{u^\gamma}.
\end{split}
\end{equation*}
  This concludes the Step.

   {\sl Step  13. $
 (\z,D \chi_{\{u>0\}})=|D \chi_{\{u>0\}}|$}

By \eqref{uuu} we need to prove the opposite inequality $(\z, D\chi_{\{u>0\}})\geq|D\chi_{\{u>0\}}|$ . We note that by Steps 11 and 12 for every
 $\varphi \in C^\infty_0 (\Omega)$ with $\varphi \ge 0$ we have
\begin{equation}\label{BVBBVV1}\nonumber
\int_\Omega \varphi |D  \chi_{\{u>0\}}| +\int_\Omega  \chi_{\{u>0\}}\, \z\cdot\nabla \varphi\, dx\leq
\int_\Omega \frac{f}{u^\gamma}\varphi \, dx=-\int_\Omega \varphi  \chi_{\{u>0\}}^*\Div\z \,.
\end{equation}
Therefore by \eqref{divergence} we get
\begin{equation}\label{BVBBVVV1}\nonumber
\int_\Omega\varphi |D \chi_{\{u>0\}}| \leq-\int_\Omega \chi_{\{u>0\}}\, \z\cdot\nabla \varphi\, dx-
\int_\Omega \chi_{\{u>0\}}^*\varphi \,\Div\z = \int_\Omega \varphi(\z,D\chi_{\{u>0\}}).
\end{equation}
The arbitrariness of $\varphi$ implies that
\begin{equation}\label{BVBBVVVB1}\nonumber
|D \chi_{\{u>0\}}| \leq (\z,D\chi_{\{u>0\}})\,,
\end{equation}
as measures on $\Omega$.

  This concludes the proof.
\end{pf}

%{{\color{red}
\begin{remark}\rm
We recall that the boundary condition used, as a rule, in the Dirichlet problem for equations involving the $1$--Laplacian is $[\z,\nu]\in \sg(-u)$ holds on $\partial\Omega $, which requires that $\z\in \DM(\Omega)$. This is not the case in our situation since we just have $\z\in \DM_{\rm{loc}}(\Omega)$ and therefore  the trace $[\z,\nu]$ is not defined and the usual boundary condition has no sense. On the contrary, the above
 boundary condition $h(u)+[h(u)\z,\nu]=0$, satisfied for every Lipschitz--continuous and nondecreasing function $h:[0,+\infty[\to [0,+\infty[$ such that $h(0)=0$,
 is completely meanigful since $h(u)\z\in \DM(\Omega)$ by Proposition \ref{poiu}.
\end{remark}

\begin{remark}\rm
We point out that if $k>0$, then $\chi_{\{u>k\}}\in BV(\Omega)$. This fact is a consequence of the inequality
\[
\frac f{u^\gamma}h(G_ku)\le \frac f{k^\gamma}h(G_ku)\in L^1(\Omega)\,,
\]
valid for every Lipschitz--continuous and nondecreasing $h:[0,+\infty[\to [0,+\infty[$ such that $h(0)=0$,
since we may apply the variational formulation to $h(G_k(s))$ and follow the argument of Step 9. Moreover, we can deduce the following facts as well:
\begin{align*}
   & -\Div(\chi_{\{u>k\}}\z)+|D\chi_{\{u>k\}}|=\frac f{u^\gamma}\chi_{\{u>k\}} \hbox{ in } \dis (\Omega)\,;\\
   &  (\z,D \chi_{\{u>k\}})=|D \chi_{\{u>k\}}|  \hbox{ as measures in } \Omega\,;\\
   &  \chi_{\{u>k\}}+[\chi_{\{u>k\}}\z,\nu]=0 \hbox{ holds } \mathcal H^{N-1}\hbox{--a.e. on } \partial\Omega\,.
\end{align*}
\end{remark}
%}}

\begin{remark}\rm
The $BV$-estimate proven in Step 2 depends only on the fact that $\gamma>0$ (see \eqref{BV1}). We point out that in the non singular setting (when $\gamma=0$) this fact does not hold, unless $f$ is small enough. This fact agrees with \cite{MST1}, where it is seen that only if $\|f\|_{W^{-1,\infty}}\leq 1$ there is a $BV$--estimate.

\end{remark}

\begin{remark}\rm
The $L^\infty$--estimate proven in Step 3 is not a consequence of Stampacchia's procedure. We remark that it depends on the singularity (since it goes to $0$ at infinity) and the fact that we are dealing with the $1$-Laplacian operator.
\end{remark}

\section{The case of strictly positive $f$} In this section we will assume that $f(x)> 0$ a.e. in $\Omega$
%{{\color{red}
and we will prove those specific features holding in this case, in particular, we will see a uniqueness result. We will see in Section 6, uniqueness does not hold for $f$ which vanishes on a set of positive Lebesgue measure.
%Nevertheless, we recall that, if $f\equiv 0$, we obtain that the trivial solution is the only one as a consequence of \cite{SWZ}.

Remark that $f\equiv 0$ must be studied in a different way, the solutions being $1$-harmonic functions. If $f\equiv 0$, then the approximate solution $u_p$ vanishes and so the solution $u$ founded in Theorem \ref{teoexist} becomes the trivial solution. Moreover, applying the results in \cite{SWZ} (which require some additional geometrical assumptions) we obtain that this trivial solution is the only one.

    \begin{Theorem}\label{teoexist1}
  Let $f\in L^N(\Omega)$ be such that $f(x)> 0$ for almost all $x\in\Omega$. Then the solution $u$  we have found in Theorem \ref{teoexist} to problem \eqref{problema1} satisfies
  \begin{enumerate}
    \item[$(a)$] $u(x)> 0$ for almost all $x\in\Omega$,
      
    \item[$(b)$] $\displaystyle \frac f{u^\gamma}\in L^1(\Omega)$,
  \end{enumerate}
   and there exists  $\z\in\DM(\Omega)$ with $\Div\z\in L^1(\Omega)$ and $\|\z\|_{\infty} \le 1$ such that
\begin{enumerate}
\item[$(c)$]  $\displaystyle -\Div \z  =\frac f{u^\gamma} \quad\text{ in }\, \dis (\Omega)\,,$

\item[$(d)$] $\displaystyle(\z, Du)=|Du| \,\text{ as measures in }\Omega\,,$

\item[$(e)$]    $[\z,\nu]\in \sg(-u)$ $\mathcal H^{N-1}$--a.e. on $\partial\Omega$.
\end{enumerate}
    \end{Theorem}

\begin{pf}
Since $u$ is a solution, we already know that (d) holds. Condition (a) is a consequence of Remark \ref{rem1}, and it yields condition (c) by Definition \ref{def1}. Observe that assuming condition (b), we deduce in a straightforward way that $\z\in\DM(\Omega)$ (with $\Div\z\in L^1(\Omega)$) and so $[\z,\nu]$ has sense and (e) holds. Hence, only condition (b) remains to be proved. We will see it in two steps.

   {\sl Step  1. For every nonnegative $v\in W_0^{1,1}(\Omega)$, we have $\int_\Omega\z\cdot\nabla v\, dx=\int_\Omega\frac f{u^\gamma}v\, dx$.}

In this Step we use an argument close to those used in Section 3, where existence and uniqueness of approximating problems \eqref{hjhjss1} are proved. We repeat these arguments for the sake of completeness (note that here $p=1$).
Consider a sequence $\varphi_n\in C^{\infty}_0(\Omega)$ such that $\varphi_n\ge0$ and
\begin{equation}\label{rrrr0}
\varphi_n\to v\quad { \hbox { strongly in } } W^{1,1}_0(\Omega).
\end{equation}

We take for every $\eta>0$ the function
\begin{equation*}\label{hjhjssdd0}
\rho_\eta*(v\wedge \varphi_n),
\end{equation*}
where $\rho_\eta$ is a standard convolution kernel and  $v\wedge \varphi_n:=\inf\{v,\varphi_n\}$. By taking it as test function in (c), we get
\begin{equation}\label{hjhjssffbb0}
    \int_\Omega \z\ \cdot \nabla (\rho_\eta*(v\wedge \varphi_n))\, dx = \int_\Omega \frac{f}{u^\gamma}\, (\rho_\eta*(v\wedge \varphi_n))\,dx\,.
\end{equation}
We want to pass to the limit as $\eta\to 0$ using
\begin{equation*}\label{rrrrss0}
\rho_\eta*(v\wedge \varphi_n)\to v\wedge \varphi_n\quad { \hbox { strongly in } } W^{1,1}(\Omega).
\end{equation*}
This implies for the left hand side of \eqref{hjhjssffbb0} that
\begin{equation}\label{hjhjssffww0}
    \int_\Omega \z \cdot \nabla (\rho_\eta*(v\wedge \varphi_n))\, dx
    \to
        \int_\Omega \z \cdot \nabla (v\wedge \varphi_n)\, dx
    \end{equation}
    As far as the right hand side of \eqref{hjhjssffbb0} is concerned, let us observe that, for $\eta>0$ small enough,
$$\sop(\rho_\eta*(v\wedge \varphi_n))\subseteq K_n,$$
where $K_n$ is a compact set contained in $\Omega$. Moreover, it follows from $\frac f{u^\gamma}\in L_{loc}^1(\Omega)$ (see \eqref{estima1}),  that
\begin{equation*}
\frac{f}{u^\gamma}\in L^1(K_n).
\end{equation*}
 Furthermore, for any $\eta>0$
\begin{equation*}
\|\rho_\eta*(v\wedge \varphi_n)\|_{L^\infty}\leq
\|v\wedge \varphi_n\|_{L^\infty}
\end{equation*}
and, as $\eta\to 0$,
\begin{equation*}
\rho_\eta*(v\wedge \varphi_n)\to
v\wedge \varphi_n\qquad a.e.,
\end{equation*}
and so
\begin{equation*}
\rho_\eta*(v\wedge \varphi_n)\rightharpoonup
v\wedge \varphi_n\qquad w^*-L^\infty.
\end{equation*}
We conclude that, as $\eta\to 0$,
    \begin{equation}\label{floflo0}
 \int_\Omega \frac{f}{u^\gamma}\, (\rho_\eta*(v\wedge \varphi_n))\,dx\to  \int_\Omega \frac{f}{u^\gamma}\, (v\wedge \varphi_n)\,dx.
\end{equation}
  By \eqref{hjhjssffbb0},   \eqref{hjhjssffww0} and \eqref{floflo0}   we obtain
      \begin{equation}\label{flafla0}
\int_\Omega \z \cdot \nabla (v\wedge \varphi_n)\, dx
=
 \int_\Omega \frac{f}{u^\gamma}\, (v\wedge \varphi_n)\,dx.
\end{equation}
    Now, we are going to pass to the limit in \eqref{flafla0}, as $n\to +\infty$. Since
    \begin{equation*}
v\wedge \varphi_n\to
v\qquad \hbox{ in }   W^{1,1}(\Omega),
\end{equation*}
we have
          \begin{equation*}\label{fluflu0}
\int_\Omega \z \cdot \nabla (v\wedge \varphi_n)\, dx
\to
\int_\Omega \z \cdot \nabla v\, dx
.
\end{equation*}
We are going now to prove that
\begin{equation*}\label{flifli0}
\int_\Omega \frac{f}{u^\gamma}\, (v\wedge \varphi_n)\,dx
\to
\int_\Omega \frac{f}{u^\gamma}\, v\,dx\,,
\end{equation*}
as $n\to +\infty$. Indeed
\begin{equation*}\label{flicflic0}
 \frac{f}{u^\gamma}\, (v\wedge \varphi_n)
\to
\frac{f}{u^\gamma}\, v\,\qquad \hbox{a.e.}
\end{equation*}
and
$$
0\leq \frac{f}{u^\gamma}\, (v\wedge \varphi_n)
\leq
\frac{f}{u^\gamma}\, v.
$$
Then by Lebesgue dominated convergence theorem, it is sufficient to prove that
\begin{equation}\label{flefle0}
\frac{f}{u^\gamma}\,v \in L^1(\Omega)\,.
\end{equation}
Indeed, by condition (c) we have
\begin{equation*}
\int_\Omega \frac{f}{u^\gamma}\varphi_n\,dx=\int_\Omega \z\cdot \nabla \varphi_n \,dx\leq \|\z\|_\infty\int_\Omega  |\nabla \varphi_n| \,dx\le C\,,
\end{equation*}
where the last inequality is due to \eqref{rrrr0}.
It follows from \eqref{rrrr0} and Fatou's Lemma that \eqref{flefle0} holds.
Therefore, Step 1 is proved.

   {\sl Step  2. The inequality $\int_\Omega\frac f{u^\gamma}v\, dx\le \int_\Omega|\nabla v|+\int_{\partial\Omega}v\,d\mathcal H^{N-1}\, dx$ holds for every nonnegative $v\in W^{1,1}(\Omega)\cap L^\infty(\Omega)$.}

Applying \cite[Lemma 5.5]{An}, we find a sequence $(w_n)_n$ in
 $W^{1, 1}(\Omega)\cap C(\Omega)$ such that
 \begin{equation*}\begin{array}{ll}
 \hbox{\rm (1) }
 w_n|_{\partial\Omega}=v|_{\partial\Omega}\,.\\
\\
 \hbox{\rm (2) } \displaystyle\int_\Omega|\nabla
 w_n|\,dx\le\displaystyle\int_{\partial\Omega}v\,d\mathcal H^{N-1}+\frac1n\,.\\
\\
 \hbox{\rm (3) } w_n(x)\to0,  \qquad\hbox{ for all } x
 \in\Omega\,.
\hskip5cm
 \end{array}\end{equation*}
 Since $|v-w_n|\in W_0^{1,1}(\Omega)$, Step 1 becomes
 \begin{multline*}
   \int_\Omega\frac f{u^\gamma}|v-w_n|\, dx=\int_\Omega\z\cdot\nabla|v-w_n|\, dx \\
   \le \|\z\|_\infty\int_\Omega|\nabla v|\, dx+\|\z\|_\infty \int_\Omega|\nabla w_n|\, dx \le \int_\Omega|\nabla v|+\int_{\partial\Omega}v\,d\mathcal H^{N-1}+\frac1n\,.
 \end{multline*}
 Then Fatou's Lemma implies
 \begin{equation*}
 \int_\Omega\frac f{u^\gamma}v\, dx\le \liminf_{n\to\infty}\int_\Omega\frac f{u^\gamma}|v-w_n|\, dx\le \int_\Omega|\nabla v|+\int_{\partial\Omega}v\,d\mathcal H^{N-1}\,.
\end{equation*}
In particular, chosing $v(x)=1$ for all $x\in\Omega$, we obtain
$\frac f{u^\gamma}\in L^1(\Omega)\,,$
as desired.
\end{pf}

As a consequence of Theorem \ref{teoexist1} (b) and (c), every $v\in BV(\Omega)\cap L^\infty(\Omega)$ can be taken as test function in our problem.

    \begin{Corollary}\label{test-bv}
  Let $f\in L^N(\Omega)$ be such that $f(x)> 0$ for almost all $x\in\Omega$. Then
  \begin{equation}\label{test}
    \int_\Omega (\z, Dv)-\int_{\partial\Omega}v[\z,\nu]\, d\mathcal H^{N-1}=\int_\Omega \frac f{u^\gamma}v\, dx
  \end{equation}
  holds for every $v\in BV(\Omega)\cap L^\infty(\Omega)$.
    \end{Corollary}

\begin{pf}
  Fix $v\in BV(\Omega)\cap L^\infty(\Omega)$.
  It follows from Theorem \ref{teoexist1} (b) and (c) that  $\hbox{div\,}\z\in L^1(\Omega)$, so that we may apply Anzellotti's theory to the case $\hbox{div\,}\z\in L^1(\Omega)$ and $v\in BV(\Omega)\cap L^\infty(\Omega)$. This gives for every  $v\in BV(\Omega)\cap L^\infty(\Omega)$
  \begin{equation*}
-\int_\Omega v\,\Div\z\, dx=\int_\Omega \frac f{u^\gamma}v\, dx,
\end{equation*}
which follows by (c) of Theorem \ref{teoexist1}. Then equation \eqref{test} is a consequence of Green's formula \eqref{GreenIII} and \eqref{des2}.
\end{pf}

    \begin{Theorem}\label{uniq}
For every $f\in L^{N}(\Omega)$, with $f(x)> 0$ a.e.,
there exists at most a solution of \eqref{problema1} in the sense of Definition \ref{def1}.
    \end{Theorem}
\begin{pf}
Assume to get a contradiction that $u_1$ and $u_2$ are positive solutions to problem \eqref{problema1}. Then, for $i=1,2$, $u_i(x)>0$ for almost all $x\in\Omega$ and there exist vector fields $\z_i$ satisfying all the requirements of Definition \ref{def1}. Applying Corollary \ref{test-bv} to the difference of the solutions and taking conditions (d) and (e) in Theorem \ref{teoexist1} into account, we obtain
\begin{gather*}
  \int_\Omega|Du_1|-\int_\Omega(\z_1, Du_2)+\int_{\partial\Omega}(u_1+u_2[\z_1,\nu])\, d\mathcal H^{N-1}=\int_\Omega\frac f{u_1^\gamma}(u_1-u_2)\, dx \\
  \int_\Omega|Du_2|-\int_\Omega(\z_2, Du_1)+\int_{\partial\Omega}(u_2+u_1[\z_2,\nu])\, d\mathcal H^{N-1}=-\int_\Omega\frac f{u_2^\gamma}(u_1-u_2)\, dx\,.
\end{gather*}
Adding both equations and rearranging its terms, we get
\begin{multline*}
  \int_\Omega\big[|Du_1|-(\z_2, Du_1)\big]+ \int_\Omega\big[|Du_2|-(\z_1, Du_2)\big]\\
  +\int_{\partial\Omega}(u_1+u_1[\z_2,\nu])\, d\mathcal H^{N-1}+ \int_{\partial\Omega}(u_2+u_2[\z_1,\nu])\, d\mathcal H^{N-1}\\
  =\int_\Omega\big(\frac f{u_1^\gamma}-\frac f{u_2^\gamma}\big)(u_1-u_2)\, dx=\int_\Omega\frac {-f}{u_1u_2}(u_1^\gamma-u_2^\gamma)(u_1-u_2)\,.
\end{multline*}
Since $\displaystyle \int_\Omega\big[|Du_1|-(\z_2, Du_1)\big]\ge0$, $\displaystyle \int_\Omega\big[|Du_2|-(\z_1, Du_2)\big]\ge0$, $\displaystyle \int_{\partial\Omega}(u_1+u_1[\z_2,\nu])\, d\mathcal H^{N-1}\ge0$ and $\displaystyle \int_{\partial\Omega}(u_2+u_2[\z_1,\nu])\, d\mathcal H^{N-1}\ge0$, we infer that
\begin{equation*}
  0\le \int_\Omega\frac {-f}{u_1u_2}(u_1^\gamma-u_2^\gamma)(u_1-u_2)\,,
\end{equation*}
whose integrand is nonpositive,
so that this term vanishes. It follows that the integrand vanishes as well. Therefore, $(u_1^\gamma-u_2^\gamma)(u_1-u_2)=0$ and so $u_1=u_2$ in $\Omega$.
\end{pf}

\section{Explicit $1$--dimensional solutions}

In this Section we will show some explicit solutions taken $\Omega=]-1,1[$ and $\gamma=1$. We restrict to the one dimensional case for the sake of simplicity.
Indeed, in any dimension and considering any $\gamma\in ]0,1]$, one may think that the constant function $u(x)=\|f\|_{W^{-1,\infty}(\Omega)}^{1/\gamma}$ is a solution to our problem reasoning as follows. Since then $f/u^\gamma=f/\|f\|_{W^{-1,\infty}(\Omega)}$, and applying  \cite[Theorem 4.3 and Remark 4.7]{MST1} we get a vector field $\z\in \DM$ satisfying  $\|\z\|_{\infty} \le 1$ and $-\Div \z  =\frac f{u^\gamma}$  in $\dis (\Omega)$. Obviously, $ (\z, Du)=0=|Du|$ holds as measures in $\Omega$. Unfortunately, the condition $ [\z,\nu]\in \sg(-u)$  $\mathcal H^{N-1}$--a.e. on $\partial\Omega$ is not guaranteed.
Nevertheless, when $\Omega=]-1,1[$, it is easy to define a function $\z$ which satisfy the weak form of boundary condition, just choosing $\z$ such that $\z(-1)=1$ if $u(-1)\ne0$ and $\z(1)=-1$ if $u(1)\ne0$.

In what follows we deal with the problem
\begin{equation}\label{1-pr}
  \left\{\begin{array}{ll}
 \displaystyle -\Big(\frac{u^\prime}{|u^\prime|}\Big)^\prime=\frac fu\,, &\hbox{in }]-1,1[\,;\\[4mm]
  u(-1)=0=u(1)\,.
  \end{array}\right.
\end{equation}
It is worth recalling the conditions of being a solution in this $1$--dimensional setting.
As just mentioned, the boundary condition becomes $\z(-1)=1$ if $u(-1)\ne0$ and $\z(1)=-1$ if $u(1)\ne0$.
On the other hand, the condition $-\z^\prime\chi_{\{u>0\}}=f/u$ implies that $\z^\prime$ is a function and so $\z$ does not jump (at least where $u>0$).
Finally, it follows from $(\z, u^\prime)=|u^\prime|$ that $\z(x)=1$ if $u$ ``increases" at $x$ and $\z(x)=-1$ if $u$ ``decreases" at $x$.

\begin{Example}\rm
Set $f(x)=\chi_{]-1/2, 1/2[}(x)$. From the previous argument, we already know that a solution to \eqref{1-pr} is given by $u_1(x)=\frac12$ with
\[
\z_1(x)=\left\{\begin{array}{ll}
1\,, &\text{if }-1<x<-\frac12\,;\\[3mm]
-2x\,,&\text{if }-\frac12\le x\le \frac12\,;\\[3mm]
-1\,,&\text{if }\frac12<x<1\,.
\end{array}\right.
\]
The same function $\z_1$ allows us to check that any $u$ satisfying
$u$ nondecreasing in $]-1,-\frac12]$, 
 $u(x)=\frac12$ in $]-\frac12,\frac12[$ and
  $u$ nonincreasing in $[\frac12,1[$,
  is a solution to \eqref{1-pr}.
For instance, $u_2(x)=\frac12\chi_{]-\frac12,\frac12[}(x)$ defines a solution.
We conclude that uniqueness does not hold.

Observe that, for $u_2=\frac12\chi_{]-\frac12,\frac12[}$, other choices of the auxiliary function are possible, namely:
\[
\z_2(x)=\left\{\begin{array}{ll}
2(x+1)\,, &\text{if }-1<x<-\frac12\,;\\[3mm]
-2x\,,&\text{if }-\frac12\le x\le \frac12\,;\\[3mm]
2(x-1)\,,&\text{if }\frac12<x<1\,.
\end{array}\right.
\]
\end{Example}

\begin{remark}\label{oss2}\rm
We point out that the function $\z$ is not unique. Moreover, what is really important in our equation is the identity
\[
-\z^\prime\chi^*_{\{u>0\}}=\frac fu\,.
\]
We remark that in the previous example $\z_1^\prime=\z_1^\prime\chi^*_{\{u>0\}}$, while $\z_2^\prime\ne\z_2^\prime\chi^*_{\{u>0\}}$.

Furthermore, one may wonder if the identity
\[
-(\z\chi_{\{u>0\}})^\prime=\frac fu
\]
holds as well. The answer is negative since one can easily check that, for solution $u=\frac12\chi_{]-\frac12,\frac12[}$, it yields
\[
(\z_1\chi_{\{u>0\}})^\prime=-2\chi_{]-\frac12,\frac12[}+\delta_{-1/2}+\delta_{1/2}=\z_1^\prime\chi^*_{\{u>0\}}+\delta_{-1/2}+\delta_{1/2}\ne \z_1^\prime\chi^*_{\{u>0\}}\,.
\]
Note that the terms $\delta_{-1/2}$ and $\delta_{1/2}$ are the measure $|D\chi_{\{u>0\}}|$ concentrated in $\partial^*\{u>0\}$.
\end{remark}

\begin{Example}\rm
Consider the function $f:]-1,1[\to\R$ defined by
\[
f(x)=\left\{\begin{array}{cl}
x-\frac12\,, &\hbox{if }\frac12\le x<1\,;\\[3mm]
0\,,&\hbox{if }-\frac12< x<\frac12\,;\\[3mm]
-x-\frac12\,,&\hbox{if }-1< x\le-\frac12\,.\\
\end{array}\right.
\]

It is straightforward that a solution of \eqref{1-pr} is given by the constant function $u_1(x)=\frac18$ with
\[
\z_1(x)=\left\{\begin{array}{cl}
-4x^2+4x-1\,, &\hbox{if }\frac12\le x<1\,;\\[3mm]
0\,,&\hbox{if }-\frac12< x<\frac12\,;\\[3mm]
4x^2+4x+1\,,&\hbox{if }-1< x\le-\frac12\,.\\
\end{array}\right.
\]

Another solution is given by $u_2(x)=\frac1{16}\chi_{\{f>0\}}(x)$ with
\[
\z_2(x)=\left\{\begin{array}{cl}
-8x^2+8x-1\,, &\hbox{if }\frac12\le x<1\,;\\[3mm]
2x\,,&\hbox{if }-\frac12< x<\frac12\,;\\[3mm]
8x^2+8x+1\,,&\hbox{if }-1< x\le-\frac12\,.\\
\end{array}\right.
\]
In this case, $-\z_2^\prime\chi_{\{u_2>0\}}=f/u_2$ holds, and $\z_2^\prime\chi_{\{u_2>0\}}\ne\z_2^\prime$. Observe that this example shows that there is not uniqueness even where
$\{f>0\}$.

There are still other solutions to \eqref{1-pr}, for instance, $u_3(x)=\frac1{16}\chi_{]-1,-\frac14[\cup]\frac14,1[}(x)$ with
\[
\z_3(x)=\left\{\begin{array}{cl}
-8x^2+8x-1\,, &\hbox{if }\frac12\le x<1\,;\\[3mm]
1\,, &\hbox{if }\frac14\le x<\frac12\,;\\[3mm]
4x\,,&\hbox{if }-\frac14< x<\frac14\,;\\[3mm]
-1\,,&\hbox{if }-\frac12< x\le-\frac14\,;\\[3mm]
8x^2+8x+1\,,&\hbox{if }-1< x\le-\frac12\,.\\
\end{array}\right.
\]

\end{Example}

\end{document}